\documentclass[fontsize=10pt]{scrreprt}

\usepackage[square,numbers,nonamebreak]{natbib}
\usepackage{setspace}

\usepackage[english]{babel}

\usepackage{chngcntr}

\usepackage[hang]{footmisc}
 at 11truept
\font\ssc=pplrc9d at 11 truept
\newcommand\qedbox{$\rlap{$\sqcap$}\sqcup$}

\usepackage[hmarginratio=1:1, bottom=2cm, top=2.5cm]{geometry}

\usepackage[automark]{scrlayer-scrpage}
\cohead{}

\let\ceheadL\cehead
\renewcommand\cehead[1]{
\ceheadL{\textnormal{#1}}
}

\cehead{Ernesto Ingrosso -- Marco Trombetti}
\lefoot{}
\rofoot{}

\usepackage{graphicx}
\usepackage[usenames,dvipsnames,svgnames,table]{xcolor}
\definecolor{Maroon}{cmyk}{0, 0.87, 0.68, 0.32}
\definecolor{RoyalBlue2}{cmyk}{80,100,0,0.1}

\renewcommand\abstract[1]{
\begin{center}
{\textbf{Abstract}}
\end{center}
{
\linespread{1.1}\fontsize{9pt}{-10pt}\selectfont #1}}

\usepackage[utf8]{inputenc}
\usepackage[sc,osf]{mathpazo}
\usepackage[euler-digits]{eulervm}
\usepackage[T1]{fontenc}
\usepackage{mathtools}

\DeclareSymbolFont{operators}{\encodingdefault}{ppl}{m}{n}
\DeclareMathAlphabet{\mathbf}{\encodingdefault}{ppl}{bx}{n}
\DeclareMathAlphabet{\mathit}{\encodingdefault}{ppl}{m}{it}

\usepackage{amsfonts}
\usepackage{amsmath}
\usepackage{ragged2e}
\usepackage{titlesec}
\usepackage{amssymb}

\renewcommand{\thesection}{\arabic{section}}
 
\titleformat{\section}{\medskip\bigskip\normalfont\Large\bf}{\thesection}{0.5em}{}
\titleformat{\subsection}{\smallskip\bigskip\normalfont\large\bf}{\thesubsection}{0.5em}{}

\usepackage{amsthm}
\newtheoremstyle{dotless}{}{}{\itshape}{}{\bfseries}{}{1em}{}

\newtheorem*{theo*}{Theorem}
\newtheorem*{rem*}{Remark}
\theoremstyle{dotless}
\newtheorem{theo}{Theorem}

\newtheorem{lem}[theo]{Lemma}

\newtheorem{cor}[theo]{Corollary}
\newtheorem{rem}[theo]{Remark}
\newtheorem{ex}[theo]{Example}
\newtheorem{quest}[theo]{Question}
\renewenvironment{proof}{\smallbreak\noindent {\sc Proof \;---\;}}{\hfill\qedbox}

\counterwithout{theo}{section}
\numberwithin{theo}{section}

\DeclareOldFontCommand{\rm}{\normalfont\rmfamily}{\mathrm}
\DeclareOldFontCommand{\sf}{\normalfont\sffamily}{\mathsf}
\DeclareOldFontCommand{\tt}{\normalfont\ttfamily}{\mathtt}
\DeclareOldFontCommand{\bf}{\normalfont\bfseries}{\mathbf}
\DeclareOldFontCommand{\it}{\normalfont\itshape}{\mathit}
\DeclareOldFontCommand{\sl}{\normalfont\slshape}{\@nomath\sl}
\DeclareOldFontCommand{\sc}{\normalfont\scshape}{\@nomath\sc}

\usepackage{cancel}

\usepackage{hyperref}

\usepackage{comment}

\usepackage{soul}
\DeclareSymbolFont{newfont}{OML}{cmm}{m}{it}%
\DeclareMathSymbol{\Varrho}{3}{newfont}{37}

\begin{document}

\centerline{\bf\LARGE Groups with dual classes of subgroups \footnote{Ingrosso and Trombetti are members of the non-profit association ‘‘AGTA --- Advances in Group Theory and Applications’’ (www.advgrouptheory.com), and are supported by GNSAGA (INdAM).}}

\medskip

\centerline{E. Ingrosso -- M. Trombetti -- A. Zitouni}

\bigskip\bigskip

\noindent\textbf{Abstract.}\quad
We completely characterize insoluble groups whose subgroups are pronormal or pronormalize every subgroup. We also completely describe $p$-groups with this property, and prove many structural results in the soluble case (for example, such groups are always of length at most $3$ and have Sylow towers).

\medskip

\noindent{\it Mathematics Subject Classification \textnormal(2020\textnormal)}: 20E07, 
20D10, 
20F16, 
20D05 

\medskip


\section{Introduction}

One of the most relevant research lines in group theory is the study of groups whose proper subgroups satisfy some property. This is particularly relevant in finite group theory, where arguments by minimal counterexample are common. In such arguments, one is naturally led to consider situations in which all {\it proper} subgroups (indeed, all strictly smaller groups) satisfy a given property (and its consequences); and so if the group does not satisfy that property (of its consequences), then it must be subject to several restrictions. Among the most notable instances of this line of research are the classifications of groups whose proper subgroups are abelian, nilpotent, supersoluble, normal and so on (see \cite{course}). For the reader’s convenience, we recall in more detail three of these examples, that concern three of the most relevant embedding properties for groups.

Groups whose (proper) subgroups are normal are called {\it Dedekind} in general, and {\it Hamiltonian} if they are non-abelian: the latter are direct products of a periodic abelian group with only odd order elements, a copy of the quaternion group $Q_8$ of order $8$, and an elementary abelian $2$-group.  The normality requirement can be weakened by asking that all subgroups of the group are permutable. Recall that a subgroup $H$ of a group $G$ is {\it permutable} if $HK=KH$ for every~\hbox{$K\leq G$.}  Groups whose (proper) subgroups are permutable are said to be {\it quasihamiltonian} and are completely characterized (see \cite{schmidt}). Another way to weaken the normality hypothesis, is by asking that all subgroups are pronormal. Recall that a subgroup $H$ of a group $G$ is pronormal if, for every $g\in G$, $H$ is conjugate to $H^g$ in $\langle H,H^g\rangle$ (in this case, we say that the element $g$ {\it pronormalizes} $H$). This concept was introduced by~Philip Hall and abstracts the common features of Sylow subgroups of finite groups, normal subgroups, and maximal subgroups. Groups whose subgroups are pronormal are described in \cite{KS}.

\medskip

Over time, classification results have become increasingly refined, addressing, for instance, groups in which {\it not all} proper subgroups satisfy a given property, but only a suitably large class of them; for example (but there are plenty of papers on the argument), in \cite{kurda}, groups whose non-abelian subgroups are normal are described. 

\medskip

In reviewing this line of research, we were led to consider seemingly natural dual conditions that appear to have been overlooked. For example, what can be said about a (finite) group in which every subgroup is normal or contained in the \textit{norm} (that is, the set of elements normalizing every subgroup)? A brief analysis shows that such a condition forces any two subgroups to permute, and hence the group is quasihamiltonian. By the known characterization of quasihamiltonian groups, it follows that these groups are direct products of $p$-groups of finite exponent of the type described in our Theorem~\ref{pgroups}. Thus, the problem admits a straightforward answer.

Turning to permutability, if one attempts to formulate an analogous dual condition, then one is led to consider groups in which every subgroup is permutable or contained in the \textit{quasi-norm}, that is, the set of elements whose cyclic subgroups permute with every subgroup. A moment's reflection shows that these are precisely the quasihamiltonian groups. In other words, permutability is sort of a ‘‘fixed point’’ for this class of problems. (Another easier ‘‘fixed point’’ is centrality.)

Thus, the only reasonable property which remains to be considered is pronormality. What about groups whose subgroups are pronormal or are contained in the {\it pro-norm}, that is, the set of all elements pronormalizing every subgroup of the group? Here, the situation is much more difficult, and it is the subject of this paper to deal with it, showing that such groups have some very nice and unexpected properties.

\medskip

Our main results with respect to the class of groups whose subgroups are pronormal or contained in the pro-norm are the following ones. Note that although most of these results apply to infinite groups, we state them below in the finite case, 
with references to their most general formulations.

\begin{itemize}
    \item A complete classification of finite insoluble groups (see Theorem \ref{locallycgraded} and the subsequent results). This is certainly our main result, and its proof is based on the structure theory we obtain in case of soluble groups.
\item A characterization of finite $p$-groups with our property (see Theorems \ref{pgroups} and \ref{PLG->LF}).
\item Every finite soluble group of our type has a Sylow tower (see Corollary \ref{generalsylowtower}), and has derived length at most $3$ (see Corollary \ref{boundforsoluble}) --- even $2$ if the group is supersoluble (see~The\-o\-rem~\ref{supersoluble}).
\item A rough structural result for finite soluble groups of our type is given in Theorem \ref{pronormform} (see Theorem also \ref{q8ipercentro}).
\item Soluble groups with our property satisfy the subnormal intersection property (see Theorem \ref{sip}). In turn, this fact has many consequences for the infinite soluble groups of our type (see the corollaries after Theorem \ref{sip}).

\end{itemize}

\section{Preliminaries}

In this section we introduce some of the basic tools we are going to exploit in our work. Most of the facts we introduce will be used in the subsequent proofs without a further notice. We also include a number of examples, with the aim of clarifying the setting and illustrating some of the underlying ideas.

\medskip

Let's start by recalling that a group is {\it locally graded} if every non-trivial finitely generated subgroup has a non-trivial finite homomorphic images. Locally graded groups are usually introduced to deal with infinite groups in a very large settings and allows us to get rid of~{\it Tarski monsters} (that is, infinite simple group whose proper non-trivial subgroups have prime order) and similar pathologies. For example, all subgroups of a Tarski monsters are pronormal, so the classification of groups with only pronormal subgroups could be an extremely difficult task in general, but in the universe of locally graded groups this can be accomplished quite well (see \cite{dGTpro}, Lemma 2.8, and \cite{KS}). It turns out that a locally graded group $G$ with only pronormal subgroups is either Dedekind or can be written as follows: $G=H\ltimes\gamma_3(G)$, where~$H$ is Dedekind, $\gamma_3(G)$ is abelian with all subgroups normal in $G$, and $\pi(H)\cap\pi(\gamma_3(G))=\emptyset$. (Actually, to get the sufficiency we should also require that the Sylow $\pi(H)$-subgroups of $G$ are complements to $\gamma_3(G)$, which is something that in finite groups always happens.)

\begin{rem}
{\rm If $\pi$ is a set of primes, and $G$ is a group, then a {\it Sylow $\pi$-subgroup} of $G$ is just a maximal $\pi$-subgroup of $G$.}
\end{rem}

Since our aim is to deal with groups in which all subgroups are pronormal or contained in the pro-norm, we should spend some words on these concepts (most of the facts we are going to cite can be found in \cite{ultimopaper}). Let's first fix some notation: the norm of a group $G$ is  denoted with $N(G)$, while the pro-norm is denoted by $N_{\operatorname{pro}}(G)$. It turns out that the pro-norm of a finite soluble group is a (characteristic) subgroup, and such is also true in polycyclic groups. This is something very useful in proving that certain types of groups are of our type. The set of all elements pronormalizing a given subgroup $H$ of a group $G$ is denoted by $P_G(H)$, and we also set $P_K(H):=P_G(H)\cap K$ for every subgroup $K$ of $G$. (We say that  subgroup $U$ {\it pronormalizes} $H$ if $U\subseteq P_G(H)$.) Even in finite soluble groups, $P_G(H)$ may not be a subgroup, and to overcome this difficulty one usually consider maximal subgroups of~$P_G(H)$ containing $H$ as a pronormal subgroup --- these are called the {\it pronormalizers} of $H$ in~$G$ (note that these may be more than one). The pro-norm, the pronormalizers, and the sets of pronormalizing elements are controlled by the finite homomorphic images in a polycyclic group.

\medskip

Another class of groups we will be needing is that of {\it critical-Sylow-tower groups} $G$, that is, of groups with no~Sy\-low tower but in which every subgroups and homomorphic image has a~Sy\-low tower. These are described in \cite{FAT}, p.365, and the finite soluble ones are of two types:

\begin{itemize}
    \item[Type 1] $|G|=p^\alpha q^\beta$ (here $p$ and $q$ are distinct primes), $G$ has a unique minimal normal subgroup~$F(G)$ and $G/F(G)$ is minimal non-nilpotent. (Here $F(G)$ is the Fitting subgroup of $G$, that is, the largest nilpotent normal subgroup of $G$.)
    
    \item[Type 2] $G=\big(\langle x_1\rangle\ltimes A_1\big)\times\ldots\times\big(\langle x_n\rangle\ltimes A_n\big)$, where $|x_i|=p_i$, $|A_i|=p_{i+1}^{r_{i+1}}$, $[x_i,A_i]\neq\{1\}$, and the~$A_i$'s are minimal normal subgroups of $G$ (here $p_1,\ldots,p_n$ are distinct primes and~\hbox{$p_{n+1}=p_1$).}
\end{itemize}

\medskip

We end this section with the promised collection of examples (and non-examples) of groups whose subgroups are pronormal or contained in the pro-norm. Since these examples can be easily deal with using GAP \cite{GAP}, we will not provide many details in every single case. Most of our notations for well-known groups is standard, but note that $\mathcal C_n$ denotes the cyclic group of order $n$, while $\mathbb Z$ the infinite cyclic group.

\begin{ex}\label{exnon}
The following groups have subgroups that are neither pronormal nor contained in the pro-norm\textnormal:
\begin{itemize}
    \item[\textnormal{(1)}] $\operatorname{Dih}(8)$. Here, $N\big(\operatorname{Dih}(8)\big)=Z\big(\operatorname{Dih}(8)\big)$. As a consequence, $\operatorname{Sym}(4)$ and $\operatorname{Dih}(\infty)$ are groups of the same type. Here, $N_{\operatorname{pro}}\big(\operatorname{Sym}(4)\big)=N_{\operatorname{pro}}\big(\operatorname{Dih}(\infty)\big)=\{1\}$. 
    \item[\textnormal{(2)}] Every minimal non-abelian group of order $p^3$ and exponent $p$ for some odd prime $p$. Here, the norm coincides with the centre.
    \item[\textnormal{(3)}] $\operatorname{Alt}(4)\times\operatorname{Alt}(4)$. Here, the pro-norm coincides with the Sylow $2$-subgroup.
    \item[\textnormal{(4)}] $\operatorname{SL}(2,3)\times\operatorname{Sym}(3)$. Here, the pro-norm is $Q_8\times\mathcal C_3$.
    \item[\textnormal{(5)}] $Q_8\times\mathcal C_4$. Here the pro-norm equals the centre.

\end{itemize}
\end{ex}
\begin{proof}
(1)\quad $\operatorname{Dih}(8)$ is not quasihamiltonian. The pro-norm $N$ of $G=\operatorname{Sym}(4)$ is obviously contained in the Wielandt subgroup of~$G$, which coincides with the Klein $4$-subgroup $V$ of~$G$. But the only element of $V$ pronormalizing $H$ is~$(1,2)(3,4)$, so $N=\{1\}$ because it is a proper~\hbox{$G$-in}\-variant subgroup of $V$.

\medskip

\noindent(2)\quad It is not quasihamiltonian, while the computation of the norm is easy.

\medskip

\noindent(3)\quad The pro-norm is easily seen to coincide with the Sylow $2$-subgroup. On the other hand, every cyclic subgroup of order $6$ is not pronormal.

\medskip

\noindent(4)\quad Every diagonal subgroup of order $3$ is neither pronormal nor even contained in the pro-norm, which is easily seen to coincide with $Q_8\times\mathcal C_3$.

\medskip

\noindent(5)\quad Write our group as $G=\langle a,b\rangle\times\langle c\rangle$, where $\langle a,b\rangle\simeq Q_8$ and $o(c)=4$. Clearly, $Z(G)\leq N_{\operatorname{pro}}(G)$. Let $x=ac$. Then $\langle x\rangle$ is not pronormal in $G$ (conjugate by $b$). On the other hand, $x$ is not contained in the pro-norm of $G$ because it does not pronormalize $\langle bc\rangle$. Now, we easily see that the pro-norm must coincide with the centre.
\end{proof}

\begin{ex}\label{exyes}
The subgroups of the following groups are pronormal or contained in the pro-norm\textnormal:
\begin{itemize}
    \item[\textnormal{(1)}] $\operatorname{Sym}(3)\times\operatorname{Sym}(3)$. Here, $|N_{\operatorname{pro}}(\operatorname{Sym}(3)\times\operatorname{Sym}(3))|=18$.
    \item[\textnormal{(2)}] $\operatorname{Alt}(5)$. Here,  $|N_{\operatorname{pro}}(\operatorname{Alt}(5))|=40$.
    \item[\textnormal{(3)}] Every minimal non-abelian group of order $p^3$ and exponent $p^2$ for some odd prime $p$. Here, the norm coincides with the set of all elements of order $p$.
    \item[\textnormal{(4)}] $\operatorname{Alt}(4)$. Here, the pro-norm coincides with the Sylow $2$-subgroup.
    \item[\textnormal{(5)}] $(E\times Q_8)\ltimes\mathcal C_3$, where the action of $E\times Q_8$ on $\mathcal C_3$ is by inversion, and $E$ is an elementary abelian $2$-group.
    \item[\textnormal{(6)}] $\operatorname{Sym}(3)\times\mathbb Z$. Here, the pro-norm is $\mathcal C_3\times\mathbb Z$.
    \item[\textnormal{(7)}] $\operatorname{Alt}(4)\times\mathbb Z$. Here, the pro-norm is $\mathcal C_2\times\mathcal C_2\times\mathbb Z$.
    \item[\textnormal{(8)}] $Q_8\ltimes (\mathcal C_3\times\mathcal C_3)$, where the action is irreducible. Here, the pro-norm is $\mathcal C_2\ltimes (\mathcal C_3\times\mathcal C_3)$ \textnormal(this also coincides with the derived subgroup\textnormal).
\end{itemize}
\end{ex}
\begin{proof}
(1)\quad Write our group as follows: $G=\big(\langle a_1\rangle\ltimes\langle b_1\rangle\big)\times\big(\langle a_2\rangle\ltimes\langle b_2\rangle\big)$, where $a_i$ has order~$2$, and $b_i$ has order $3$ for every $i\in\{1,2\}$. Clearly, $a_i$ does not pronormalize $\langle b_1b_2\rangle$, while $b_1,b_2$ and $a_1a_2$ are contained in the pronorm. Thus, a subgroup $H$ that is not contained in the pro-norm must contain a conjugate of $a_1$ or $a_2$. Suppose $H$ contains $a_1$. If it contains $b_1$, then we can factor $\langle a_1,b_1\rangle$ out and note that all subgroups of $\operatorname{Sym}(3)$ are pronormal. If it does not contain $b_1$, but contains $b_2$, then we factor by $\langle b_2\rangle$ and note that again all subgroups of the quotient are pronormal. If $H$ does not contain neither $b_1$ nor $b_2$, then we may assume $H\leq\langle a_1\rangle\times \langle a_2\rangle$, and that actually either $H=\langle a_i\rangle$ for some $i$, or $H=\langle a_1,a_2\rangle$. In both cases, $H$ is trivially pronormal in $G$.

\medskip

\noindent(2)\quad Use GAP (or see Theorem \ref{5mproof}).

\medskip

\noindent(3)--(4)\quad Easy.

\medskip

\noindent(5)\quad All subgroups are pronormal in this case (see the main theorem of \cite{KS}).

\medskip

\noindent(6)\quad Our group $G$ is polycyclic, so the pro-norm is a subgroup of $G$. Write $G=(\langle a\rangle\ltimes \langle b\rangle)\times \langle c\rangle$, where $o(a)=2$, $o(b)=3$, and $o(c)=\infty$. In order to prove that $N_{\operatorname{pro}}(G)=\langle b,c\rangle$, we only need to show that $b$ is contained in the pro-norm ($c$ being central, and $G$ not having only pronormal subgroups). Clearly, $b$ centralizes $\langle b,c\rangle$. Let $H$ be any subgroup of $G$ that is not contained in~\hbox{$\langle b,c\rangle$}, and let $g$ be an element of $H$ that is not contained in~\hbox{$\langle b,c\rangle$.} Then $g^2\in \langle c\rangle$, and in order to show that $b$ pronormalizes $H$, we may assume that $c^2=g^2=1$, and so in particular that $G$ is finite and $|H|\leq 4$. But now it is very easy to realize that $b$ pronormalizes~$H$.

Finally, if $H$ is any subgroup of $G$ which is not contained in $N_{\operatorname{pro}}(G)$, then $H$ is easily seen to be pronormal in $G$ (essentially because one can assume that $o(c)=2$, which in turns implies that all subgroups of $G$ are pronormal).

\medskip

\noindent(7)\quad Again, our group $G$ is polycyclic, so the pro-norm is a subgroup of $G$. Write $$G=\langle a\rangle\ltimes\big(\langle b\rangle\times\langle c\rangle\big)\times\langle d\rangle,$$ where $o(a)=3$, $o(b)=o(c)=2$, and $o(d)=\infty$. As in (6), to prove that $N_{\operatorname{pro}}(G)=\langle b,c,d\rangle$, we may assume $o(d)=3$, which immediately implies the claim (see also Lemma \ref{pnormldedek}).

Finally, let $H$ be any subgroup of $G$ which is not contained in $N_{\operatorname{pro}}(G)$. In order to show that~$H$ is pronormal, we may assume $o(d)=3$. Clearly, if $H$ contains a $3$-element, then $H\trianglelefteq G$. If not then $H$ is a $2$-subgroup, so either $|H|=2$ or $|H|=4$. In both cases, $H$ is clearly pronormal in $G$.

\medskip

\noindent(8)\quad This is an easy exercise (see also Lemma \ref{pnormldedek}).~\end{proof}

\section{The finite soluble case}

This section is devoted to proving general structural results in the soluble setting for finite groups whose subgroups are pronormal or contained in the pro-norm. The first result provides a complete description of the locally finite $p$-groups of this type for any prime $p$. This, together with~Co\-rol\-la\-ry~\ref{cortorsionfree}, yields for example that any (locally) nilpotent group of our type has nilpotency class at most $2$.

\begin{theo}\label{pgroups}
Let $p$ be a prime, and $G$ a non-abelian locally finite $p$-group. Then all subgroups of $G$ are normal or contained in the norm if and only if $G$ is one of the following groups:
\begin{itemize}
    \item[\textnormal{(1)}] $Q_8\times E$, where $E$ is an elementary abelian $2$-group.
    \item[\textnormal{(2)}] $G=\langle x\rangle\ltimes\big(\langle a_1\rangle\times B\big)$, where $G'\leq\langle a_1\rangle$, $B\leq Z(G)$, and $o(a_1)\geq |G'|^2/p\cdot \operatorname{max}\{o(x),\operatorname{exp}(B)\}$. In particular, $G'\leq Z(G)$.
\end{itemize}
\end{theo}
\begin{proof}
Suppose first that all subgroups of $G$ are normal or contained in the norm. Every finite subgroup of $G$ is quasihamiltonian, so $G$ is quasihamiltonian too.  Thus, if $Q_8$ is involved in $G$, then $G$ is of type (1). Assume $Q_8$ is not involved in $G$. Then there is a finite-exponent abelian subgroup~$A$ of $G$ and an element $x\in G$ such that $G=\langle x\rangle A$ with $x$ acting on $A$ as a power automorphism which is~\hbox{$\equiv_41$} in case $p=2$ (see Theorem 2.3.1 of \cite{schmidt}).

\smallskip

From now on, since the proof is exactly the same and for reader's convenience, we prove the statement assuming that $G$ is finite.

\smallskip

Write $A=\langle a_1\rangle\times\langle a_2\rangle\times\ldots\times\langle a_n\rangle$, where $o(a_1)\geq o(a_2)\geq\ldots \geq o(a_n)$. We claim that $\langle a_2,\ldots,a_n\rangle\leq Z(G)$.  Suppose not. Clearly, we may assume $n=2$. If $o(a_1)>o(a_2)$, then by induction on the order, $\langle a_2\rangle\big(\langle a_1a_2\rangle[p]\big)/\big(\langle a_1a_2\rangle[p]\big)$ and $\langle a_2\rangle/\big(\langle a_2\rangle[p]\big)$ are centralized by~$x$, so $a_2$ is centralized by $x$. Thus, $o(a_1)=o(a_2)$. Assume that $o(x) > o(a_1)$. Since $\langle x \rangle \cap A$ is a cyclic subgroup of~$A$, we may assume without loss of generality that $\langle x \rangle \cap A \leq \langle a_1 \rangle$. It follows that $\langle xa_2 \rangle$ is not normalized by $a_1$. However, $xa_2$ does not belong to the norm because it fails to normalize $\langle x \rangle$, which yields a contradiction.
Thus, $o(x)\leq o(a_1)$, and 
we may even assume that $\langle x\rangle\cap A=\{1\}$. Now, $\langle a_1x\rangle$ is not normal in $G$ because it is not normalized by $a_2$, so it must be contained in the norm of~$G$. However,~$a_1x$ does not normalize $\langle x\rangle$, so we obtain a contradiction. The claim is proved. In particular, $o(a_1)\geq o(a_2)|G'|$.

We claim that $x$ and $a_1$ can be chosen in such a way that $o(x)<o(a_1)$. Let $L=\langle x,a_1\rangle$. If $|L'|=p$, then $L$ is minimal non-abelian and metacyclic, so we may write $L=\langle x'\rangle\ltimes\langle a_1'\rangle$ for suitable elements $x',a_1'$. In this case, $o(x')<o(a_1')$ otherwise $\langle x'a_1'\rangle$ is not normal and does not even normalize $\langle x'\rangle$. Similarly,
if $M=\langle x\rangle\cap\langle a_1\rangle=\{1\}$, then $o(x)<o(a_1)$ since otherwise~$\langle xa_1\rangle$ is not normal and does not even normalize $\langle x\rangle$ (put $x'=x$ and $a_1'=a_1$ in this case). Thus, we may assume $|L'|\geq p^2$ and $M\neq\{1\}$. Let $S$ be the socle of $M$. By induction we may find $x'$ and $a_1'$ such that $L/S=\langle x'S\rangle\ltimes\langle a_1'S\rangle$ and $o(x'S)<o(a_1'S)$. If $S\not\leq\langle a_1'\rangle$, then $G'$ is not cyclic, a contradiction. Thus, $S\leq\langle a_1'\rangle$, so $\langle a_1'\rangle$ is normalized by $x'$ and $o(a_1')>o(x')$. In any case, $o(a_1')\geq o(a_1)$.

Now, $A_1=\langle a_1',a_2,\ldots,a_n\rangle$ is abelian and $x'$ acts on $A_1$ as a power automorphism because $o(a_1')\geq o(a_1)\geq o(a_2)|G'|$. Finally, by suitably replacing $x'$, we may further require that $\langle x'\rangle\cap A_1=\{1\}$ (it is in fact enough choosing $x'$ of smallest possible order).

\medskip

It remains to prove the bound in point~(2). Set $o(x) = p^t$,  $o(a_1) = p^n$, and $a_1^x = a_1^{1+p^s}$ for some $s\ge 1$. In particular, $\langle a_1^{p^{n-s}}\rangle= Z(G)$ and $|\langle a_1 \rangle : Z(G)| = |G'|$. Assume $o(x) > |Z(G)|$. If $y$ is an element of $\langle a_1\rangle$ with $o(y)=o(x)$, then $\langle xy\rangle$ is not normalized by $a_1$ and does not normalize $\langle x\rangle$, a contradiction. Therefore, $o(x)\le |Z(G)|=o(a_1)/|G'|$.  Now, assume that $o(x) > o(a_1)p/|G'|^2$. Let $y=a_1^{p^{(n-t)-(n-s)}}$. Note that $o(a^{p^{n-t}})=p^t$ and that $\langle y\rangle>Z(G)$. It follows that $\langle xy \rangle$ does not normalize $\langle x \rangle$ and is not normalized by $a_1$ because the assumption yields that $\langle xy\rangle\not\geq G'$. 
If we replace $x$ by $xa_i$, $i>1$, in this argument, then we see that all the required conditions are satisfied.




\medskip

Conversely, suppose $G$ is a group of type (2). Then the conditions imply that every cyclic subgroup of $G$ which is not centralized by $x$ must contain $G'$. Thus, $x\in N(G)$. On the other hand, every cyclic subgroup which lies outside $\langle x\rangle\big(\langle a_1,\ldots,a_n\rangle\cap Z(G)\big)\leq N(G)$ contains $G'$, so $G$ satisfies the required condition.
\end{proof}

\medskip

The relevance of the previous characterization in the general theory of the class of groups we are studying can immediately be seen from the following result.

\begin{theo}\label{semi=diretto}
Let $G$ be a finite group whose subgroups are pronormal or contained in the pro-norm. If~\hbox{$G=H\ltimes P$,} where $P$ is a non-Dedekind $p$-group for some prime $p$, and $(|H|,p)=1$, then $G=H\times P$.
\end{theo}
\begin{proof}
Since $P$ is not Dedekind, it must satisfy condition (2) in the statement of Theorem~\ref{pgroups}. Write $P=\langle x\rangle\ltimes\big(\langle a_1\rangle\times\ldots\times\langle a_n\rangle\big)$, where $o(a_2)\geq\ldots \geq o(a_n)$, $P'\leq\langle a_1\rangle$, $\langle a_2,\ldots,a_n\rangle\leq Z(P)$, and $o(a_1)\geq |P'|^2/p\cdot \operatorname{max}\{o(x),o(a_2)\}$. Clearly, $a_1$ does not normalize~$\langle x\rangle$, so $\langle a_1\rangle$ is pronormal in $G$, and hence (being subnormal) also normal. Similarly, $\langle xa_1\rangle$ and $\langle a_1a_i\rangle$, $i>1$, are normal in $G$. Moreover, if $H$ centralizes $\langle a_1\rangle$ it must also centralize $\langle xa_1\rangle$ and $\langle a_1a_i\rangle$, $i>1$, because the socles are the same, and $(|H|,p)=1$, and hence $H$ centralizes $P$. Thus, we may assume that $C_H(a_1)=\{1\}<H$ and that $p$ is odd (because the automorphism group of a cyclic $2$-group is a $2$-group).

If $|P'|>p$, then induction yields that  $[H,P]\leq\Omega(P')$, so $H$ centralizes $\langle a_1\rangle$, a contradiction. Thus, $|P'|=p$.  

Now, $U=\Omega(\langle x,a_1\rangle)$ is elementary abelian of order $p^2$, so $U=U_1\ltimes \langle a_1\rangle$, where $U_1$ is $H$-invariant. If $|U_1|\neq|\langle x\rangle|$, then $U_1\trianglelefteq P$, so $U_1\trianglelefteq G$, and hence $[H,P]\leq U_1$ by induction; but this yields again the contradiction that $H$ centralizes $a_1$. Thus, we may assume $U_1=\langle x\rangle$. Similarly, we can get rid of the $a_i$ with $i\geq2$, and so assume $n=1$.

Next, observe that every subgroup of $\langle x,a_1\rangle/\langle a_1^p\rangle$ is $H$-invariant, so $H$ acts as a group of power automorphisms on it. But $\langle a_1\rangle=C_G(a_1)$, so $G/\langle a_1\rangle$ is abelian, and hence $[H,U_1]=\{1\}$. Thus, $H$ must centralize $a_1$, the final contradiction.
\end{proof}

\medskip

Using the previous result, we can now prove three main general structure results (see Corollaries \ref{generalsylowtower} and \ref{sollength}, and Theorem \ref{supersoluble}).

\begin{theo}\label{sylowtowerth}
Let $G$ be a finite soluble group whose subgroups are pronormal or contained in the pro-norm. Then $G$ has a Sylow tower.
\end{theo}
\begin{proof}
It is enough to show that the critical-Sylow-tower groups $G$ do not satisfy the hypotheses. Suppose first that $G$ is critical-Sylow-tower of type $2$. Then there are primes $p_1,\ldots,p_n,p_{n+1}=p_1$ such that $$G=\big(\langle x_1\rangle\ltimes A_1\big)\times\ldots\times\big(\langle x_n\rangle\ltimes A_n\big),$$ $|x_i|=p_i$, $|A_i|=p_{i+1}^{r_{i+1}}$, $[x_i,A_i]\neq\{1\}$, and the $A_i$ are minimal normal subgroups of $G$. Clearly, not every $A_i$ can be cyclic, otherwise we would have $p_1<p_2<\ldots<p_n<p_1$. Let $j$ be such that $A_j$ is not cyclic. Let $a\in A_n\setminus\{1\}$ and define $g$ as $x_1x_ja$ if $j\neq1$, or as $x_1a$ if $j=1$.

Put $H=\langle g\rangle$. Clearly, $H$ is not pronormal in $G$ since it is subnormal and non-normal in~\hbox{$\langle x_2,x_1,a\rangle$.} But $H$ is not even contained in the pro-norm because otherwise $A_j$ would have been cyclic.

\medskip

Assume $G$ is critical-Sylow-tower of type $1$. Thus, $|G|=p^\alpha q^\beta$ (here $p$ and $q$ are distinct primes), $G$ has a unique minimal normal subgroup~$F(G)$ and $G/F(G)$ is minimal non-nilpotent. Suppose $F(G)$ is an (necessarily elementary abelian) $p$-group, and let $P$ and $Q$ be respectively a~Sy\-low~\hbox{$p$-sub}\-group and a Sylow $q$-subgroup of $G$; in particular, $G/F(G)=P/F(G)\ltimes QF(G)/F(G)$ and $P/F(G)$ is cyclic.

If $P$ is Dedekind, then $F(G)\leq Z(P)$, so $C_G(F(G))=G$, and hence $Q$ is characteristic in~$QF(G)$, a contradiction. Thus, $P$ is not Dedekind, and so $P$ satisfies the structure of Theorem \ref{semi=diretto}.

Since $F(G)$ is a minimal normal subgroup of $G$ and $Z(P)\neq F(G)$, we have that~\hbox{$Z(P)\cap F(G)$} cannot be normalized by $Q$. Now, let $a\in Z(P)\cap F(G)$ and $x\in Q$. Then $a^x=a'$ for some $a'\in F(G)$. Set $P/F(G)=\langle yF(G)\rangle$. If we now apply conjugation by $y$, we find that $a^{x^y}=(a')^y$, where $x^yF(G)=x^\ell F(G)$ for some positive integer $1<\ell<q$. Thus, $a^{x^{y}}=a^x$ modulo~$P'$ and so $a^{x^{y-1}}\in P'$ for every $a\in Z(P)$. But $P'\leq Z(P)$ and $[Q,\langle y\rangle]F(G)=QF(G)$, so $Z(P)$ is normal in $G$. Therefore~\hbox{$Z(P)=P$,} the final contradiction.
\end{proof}

\begin{cor}\label{generalsylowtower}
Let $G$ be a locally finite soluble group whose subgroups are  pronormal or contained in the pro-norm. Then $G$ has a finite normal series with nilpotent coprime factors. In particular, $G$ has an ascending Sylow tower.
\end{cor}
\begin{proof}
A standard inverse limit argument (using Theorem \ref{sylowtowerth}) shows that $G$ has a normal series with coprime factors. If $H/K$ is a factor of this series such that $H$ is not contained in the pro-norm of $G$, then $G/H$ has only pronormal subgroups; in particular, $$\pi(G/\gamma_3(G)H)\cap\pi(\gamma_3(G)H/H)=\emptyset.$$ Let $W$ be the intersection of all the factors of the previous type. Clearly, $\pi(W)\cap\pi(G/W)=\emptyset$ and \hbox{$W\leq N_{\operatorname{pro}}(G)$.} But also $\pi(\gamma_3(G)W/W)\cap\pi(G/\gamma_3(G)W)=\emptyset$, and so we are done.
\end{proof}

\begin{cor}\label{sollength}
Let $G$ be a finite soluble group whose subgroups are pronormal or contained in the pro-norm. Then $G$ has derived length at most $3$.
\end{cor}
\begin{proof}
By Theorem \ref{sylowtowerth}, $G$ has a Sylow tower. Let $P$ be a normal Sylow $p$-subgroup of $G$ for some prime $p$. Clearly, we may assume by induction that $G$ has only one minimal normal subgroup $M$, which is necessarily a $p$-group.

Let $Q$ be a Sylow $p'$-subgroup of $G$, so $G=Q\ltimes P$. Let $X$ be any subgroup of $Q$ that is not pronormal in $Q$ (and so in $G$). Then $X$ is contained in the pro-norm $N$ of $G$. But $Z(N)N'\leq P$, so $N/(P\cap N)$ is cyclic, and $X(P\cap N)/(P\cap N)$ is normal in $G/(P\cap N)$. Then $XP/P\trianglelefteq G/P$ and so $X\trianglelefteq Q$, a contradiction. Thus, $Q$ is metabelian.

Now, suppose that $P$ is non-abelian. Then Theorem \ref{semi=diretto} shows that $P$ is a Dedekind $2$-group. By Maschke's Theorem, the socle of $P$ is cyclic, so $P$ is a quaternion group of order $8$. Then $C_G(P)=Z(P)$ because of the uniqueness of $M$, and hence $G/P$ is a proper subgroup of $\operatorname{Sym}(4)$ (which does not satisfy our hypotheses), so $G/P$ is metabelian and $G^{(3)}=\{1\}$.
\end{proof}

\begin{theo}\label{supersoluble}
Let $G$ be a finite group whose subgroups are pronormal or contained in the pro-norm. If $G$ is supersoluble, then $G$ is metabelian.
\end{theo}
\begin{proof}
Let $G$ be a minimal counterexample. Then $G$ has only one minimal normal subgroup $N$, which has prime order $p$. Clearly, $p>2$ and $G$ is not a $p$-group. Let $P=O_p(G)$ and let $C$ be the intersection of all the centralizers of the chief factors of $G$ contained in $P$. Then $P\leq C$, and $G/C$ is abelian of order prime with $p$. If $C\neq P$, then there is a minimal normal subgroup $M/P$ of $G/P$ contained in $C/P$ which is a $q$-group for some prime $q\neq p$, so~\hbox{$M=P\times L$,} where~$L$ is a non-trivial normal $q$-subgroup of $G$, a contradiction. Thus, $C=P$, and so $P$ is non-abelian. But then Theorem \ref{semi=diretto} shows that $G=Q\times P$, so $G$ is metabelian, the final contradiction.
\end{proof}

\medskip

Giving a very precise description of finite soluble groups whose subgroups are pronormal or contained in the pro-norm is not an easy task, due to the large range of possibilities that may arise. In the remainder of this section, we try at least to give a very approximate description of such groups. We start with some useful lemmas that already shed some light on the structure.

\begin{lem}\label{lempronormalizers}
Let $G$ be a group, and $H$, $K$ subgroups of $G$ with $N_G(K) \ge H$. If $g \in P_G(H) \cap N_G(K)$, then $g \in P_G(HK)$.
\end{lem}
\begin{proof}
Since $g \in P_G(H)\cap N_G(K)$, there exists $u \in \langle H, H^g \rangle\leq N_G(K)$ such that $H^u = H^g$. Thus, $gu^{-1} \in N_G(K)$, so $(HK)^{gu^{-1}} = HK$, and hence $(HK)^g = (HK)^u$. Therefore $g\in P_G(HK)$, and the statement is proved.
\end{proof}

\begin{lem}\label{pnormldedek}
Let $\pi$ be a set of primes, and $G$ a finite soluble group. If the Sylow $\pi$-subgroups of $G$ are~De\-de\-kind, then $O_\pi(G)\leq \operatorname{N}_{\mathrm{pro}}(G)$.
\end{lem}
\begin{proof}
Let $H$ be a subgroup of $G$. Then $H = H_{\pi'}H_\pi$, where $H_\pi'$ is a Sylow $\pi'$-subgroup of~$H$, and $H_\pi$ is a Sylow $\pi$-subgroup of $H$. Now, $G$ is soluble, so $O_\pi(G)$ pronormalizes $H_{\pi'}$. Consequently, if $x\in O_\pi(G)$, then there exists $u \in \langle H, H^x\rangle \le HO_\pi(G)$ such that $H_{\pi'}^x = H_{\pi'}^u$; clearly, we may assume $u\in O_\pi(G)$. Then $H^x = H_{\pi'}^xH_\pi^x =H_{\pi'}^xH_\pi =H_{\pi'}^uH_\pi^u = H^u$ because~$O_\pi(G)$ normalizes $H_\pi$. The statement is proved.
\end{proof}

\begin{lem}\label{q8lem}
Let $G=H\ltimes K$ be a finite group, where $K$ is a non-abelian Dedekind $2$-group, $H$ is abelian, and $\pi(H)\cap\pi(K)=\emptyset$. If $G$ is metabelian, then $H$ centralizes a copy of $Q_8$ contained in $K$.
\end{lem}
\begin{proof}
Suppose first that $H=\langle h\rangle$ is cyclic. Write $K=Q\times W$, where $Q=\langle a,b\rangle$ is a copy of the quaternion group of order~$8$, and~$W$ is elementary abelian. Set $Z=Z(Q)\times W$. Since $K/Z$ is a Klein $4$-group, and $G'$ is abelian, we have that $[H,K]\leq Z$. Write $Z=C_K(H)\times [H,Z]$, and put $a^h=aa_1a_2$ and $b^h=bb_1b_2$, where $a_1,b_1\in C_K(H)$ and $a_2,b_2\in [H,Z]$. Replacing $h$ by~$h^2$, we have that $a_1=b_1=1$. Thus, $a^h=aa_2$ and $b^h=bb_2$. 

Now, note that $[H,Z]=\{[h,z]\,:\, z\in Z\}$, so there exist $z_1,z_2\in Z$ such that $[h,z_1]=a_2$ and $[h,z_2]=b_2$. Let $Q_h=\langle az_1,bz_2\rangle$ is centralized by $H$; in fact, for example, $$(az_1)^h=aa_2z_1^h=aa_2z_1[h,z_1]^{-1}=az_1.$$ This proves the statement when $H$ is cyclic.

\medskip

Finally, we turn to the general case. Write $H=H_1\times\langle h\rangle$. Then $h$ centralizes not only $Q_1$ but also $Q_1^H$. Let $G_1=H_1\ltimes Q_1^H$. By induction, we see that $Q_1^H$ contains a copy $Q_2$ of $Q_8$ that is centralized by $H_1$. Consequently, $Q_2$ is centralized by $H$ and the statement is proved.
\end{proof}

\medskip

In the following if $G$ is a group, and $N$ and $H$ a normal subgroup and a subgroup of $G$, respectively, then $\operatorname{PAut}_H(N)$ denotes the subgroup of all elements of $H$ that act as power automorphisms on $N$.

\begin{lem}\label{semidiretto}
Let $G = Q\ltimes V$ be a finite soluble group, where $V$ is  abelian, and $\pi(V) \cap \pi(Q) = \emptyset$. If every subgroup of $G$ is pronormal or contained in the pro-norm then the following conditions hold\textnormal:
\begin{itemize}
        \item[\textnormal{(1)}] every non-pronormal subgroup $K$ of $Q$ is contained in the pro-norm of $G$\textnormal; 
        \item[\textnormal{(2)}] if $K$ is a pronormal subgroup of $Q$ which is not contained in $\operatorname{N}_{\operatorname{pro}}(G)$, then every subgroup $U$ of~$V$ normalized by $K$ is normalized by~$N_G(K)$.
    \item[\textnormal{(3)}] $Q/\operatorname{PAut}_Q(V)$ has only pronormal subgroups, so $Q''\leq \operatorname{PAut}_Q(V)$.
    \item[\textnormal{(4)}] $[\operatorname{PAut}_Q(V),Q]\leq C_Q(V)$.
    \item[\textnormal{(5)}] The Sylow subgroups of $Q/\operatorname{PAut}_Q(V)$ are abelian.
    \item[\textnormal{(6)}] If $V$ is an irreducible $Q$-module, and $Q$ is a $q$-group for some prime $q$, then either $Q/C_Q(V)$ is cyclic, or $Q$ is~De\-de\-kind and $Q/C_Q(V)\simeq Q_8$.
\end{itemize}

\noindent Conversely, if $G$ satisfies \textnormal{(1)} and \textnormal{(2)}, then all its subgroups are pronormal or contained in the pro-norm.
\end{lem}
\begin{proof}
Assume first that every subgroup of $G$ is pronormal in $G$ or contained in the pro-norm of $G$.  Obviously, if $K$ is non-pronormal in $Q$, then $K$ is non-pronormal in $G$. Consequently, 
   \hbox{$K \le \operatorname{N}_{\operatorname{pro}}(G)$}. 
Now, assume that $K \not \le \operatorname{N}_{\operatorname{pro}}(G)$.
   Let $U$ be a subgroup of $V$ which is~\hbox{$K$-in}\-va\-riant. Therefore, $K \ltimes U$ is pronormal in $G$ (because $K$ is not contained in the pro-norm of~$G$). 
   If~\hbox{$g\in N_G(K)$}, then there exists $u \in \langle KU, KU^g \rangle = K \ltimes \langle U, U^g \rangle$ such that $(KU)^u = KU^g$. Write $u=xy$, where $x\in K$ and $y\in\langle U,U^g\rangle$. Then $KU^g=(KU)^u=(KU)^y=K^yU$ and hence $U^g=U$, as desired. Thus, (1) and (2) are proved.

\medskip

Set $P=\operatorname{PAut}_Q(V)$. The fact that $Q/P$ has only pronormal subgroups follows from the observation that any non-pronormal subgroup of $Q$ must be contained in the pro-norm of $G$ and so it must act as power automorphisms on $\gamma_3(X)$. This proves (3). Moreover, $P/C_Q(V)$ is clearly central in $Q/C_Q(V)$, so $[P,Q]\leq C_Q(V)$. This proves (4).

\smallskip

In order to prove (5), suppose that $Q/P$ contains a copy  of the quaternion group of order $8$. Since the Sylow subgroups of $Q$ are either Dedekind or of type described in Theorem \ref{pgroups},  there exists a subgroup $W$ of $Q$ such that $W\simeq Q_8$ and $W\cap P=\{1\}$. The latter condition implies that $Z(W)$ is not contained in the pro-norm of $G$, so $W$ normalizes every subgroup of $V$ which is normalized by $Z(W)$. In turns, this shows that every cyclic subgroup of $V$ normalized by~$W$ must be centralized by $Z(W)$, and so that $Z(W)$ centralizes $V$, a contradiction. Therefore the~Sy\-low subgroups of $Q/P$ are abelian. This proves (5).

\smallskip

Finally, (6) follows at once from Theorem \ref{pgroups} and the fact that the centre of an irreducible subgroup of $\operatorname{GL}(n,p)$ is cyclic.


\medskip

Conversely, assume that $G$ satisfies (1) and (2). Since $\pi(V) \cap \pi(Q) = \emptyset$ and $V$ is abelian, so $V \le N_{\operatorname{pro}}(G)$ by~Lem\-ma~\ref{pnormldedek}.  Let $M$ be the subgroup generated by the non-pronormal subgroups of $Q$. By (1), $M\leq N_{\operatorname{pro}}(G)$.
 
Now, let $H$ be a subgroup of $G$ which is not contained in ${N}_{\operatorname{pro}}(G)$. Write $H = H_Q \ltimes U$, where $U \le V$, and $H_Q$ is the Sylow $\pi(Q)$-subgroup of $H$. We claim that $H$ is pronormal in~$G$. Replacing $H$ with a conjugate (if necessary), we may assume that $H_Q \le Q$. 
    If $H_Q$ is not pronormal in $Q$, then $H_Q \le M$, so $H \le {N}_{\operatorname{pro}}(G)$, which is a contradiction. Thus, $H_Q$ is pronormal in $Q$. On the other hand, $H_Q $ is pronormal in $H_Q \ltimes V$, since $V$ contained in the pro-norm of $G$. Therefore, $H_Q$ is pronormal in $G=Q(H_QV)$ (see Theorem 3.45 of \cite{ultimopaper}). Finally, if
 $g \in G$, then there exists $u \in \langle H_Q, H_Q^g \rangle \le \langle H, H^g \rangle$ such that $H_Q^u = H_Q^g$; but $gu^{-1} \in N_G(H_Q)$, and so $H^g = H^u$ by (2). This shows that $H$ is pronormal in $G$, and we are done. 
\end{proof}

\begin{rem}\label{remsemidiretto}
{\rm In the context of point (2) of Lemma \ref{semidiretto}, we note that if $Q$ is nilpotent, then every subgroup of $Q$ is  normal in $Q$ or contained in the pro-norm of $G$. Thus, if $K$ is a subgroup of $Q$ that does not normalize a cyclic subgroup of $V$, then $K$ and $Q$ normalize the same subgroups of $V$, and in particular  $V$ is $Q$-irreducible if and only if it is $K$-irreducible.}
\end{rem}

\begin{theo}\label{pronormform}
Let $G$ be a finite soluble group whose subgroups are  pronormal or contained in the pro-norm. Assume that not every subgroup of $G$ is pronormal. Then $G = (T \ltimes S)\ltimes X$, where $S$ is a~Sy\-low~\hbox{$p$-sub}\-group of $G$ \textnormal(here $p$ is a prime\textnormal) that is not contained in the pro-norm of $G$, all subgroups of~$T$ are pronormal in $G$, $X\leq N_{\operatorname{pro}}(G)$, $\pi(X)\cap\pi(T)=\emptyset$, and $p\not\in\pi(X)\cup\pi(T)$. Clearly, $S$ is either~De\-de\-kind, or of the type described in \textnormal{Theorem \ref{pgroups}}. 
Write $X=M\ltimes \gamma_3(X)$ and replace $S$ and $T$ in such a way that  $G=W\ltimes \gamma_3(X)$, where $W=(T\ltimes S)\ltimes M$. Then $[M,W]\leq C_M(\gamma_3(X))$ and $W''\leq\operatorname{PAut}_W(\gamma_3(X))$.

\begin{itemize}
    \item In particular, if $\pi([M,W])\cap\pi(C_M(W))=\emptyset$, then we may assume $X$ to be Dedekind.
\end{itemize}

\medskip

\noindent Assume $G$ is metabelian. Then $T$ can be chosen to be Dedekind, $G'\leq N_{\operatorname{pro}}(G)$, and $[T,S]=\{1\}$. 

\begin{itemize}
    \item[\textnormal{($\star$)}] Assume $[S,M]=\{1\}$. Then\textnormal:
    \begin{itemize}
        \item $G=S\ltimes (TX)$,
        \item all subgroups of $TX$ are pronormal, and in fact $TX=R\ltimes\gamma_3(TX)$, where \hbox{$R=T\times C_M(T)$} and $\gamma_3(TX)=[M,T]\times\gamma_3(X)$ \textnormal(in particular, $[S,R]=\{1\}$, while $\pi(T)$, $\pi(C_M(T))$, $\pi([M,T])$, and $\pi(\gamma_3(X))$ have pairwise empty intersection\textnormal),
        \item if~$S_1$ is a subgroup of~$S$ that is not contained in $\operatorname{N}_{\operatorname{pro}}(G)$, then $N_S(S_1)$ normalizes every subgroup of $\gamma_3(X)$ that is normalized by $S_1$,
        \item the non-pronormal subgroups of $S$ act as power automorphisms on $\gamma_3(X)$.
    \end{itemize} 
\end{itemize}

Conversely, if $G$ has the structure described in \textnormal{($\star$)}, then every subgroup of $G$ is pronormal or contained in the pro-norm.

\end{theo}
\begin{proof}
By Theorem \ref{sylowtowerth}, $G$ has a Sylow tower  $\{1\} =P_0 \trianglelefteq P_1 \trianglelefteq \ldots \trianglelefteq P_t = G,$ where $P_{i+1}/P_{i}$ is a Sylow $p_i$-subgroup of $G/P_i$ (here $p_i$ is a prime) for $i\in\{0, \ldots, t-1\}$. Let~$X$ be the largest term of this series that is contained in $N_{\operatorname{pro}}(G)$. Then we can write $G=T\ltimes (S\ltimes X)$, for some subgroup $T$ of $G$, and for some Sylow $p$-subgroup $S$ of $G$ (here $p$ is a prime) that is not contained in $N_{\operatorname{pro}}(G)$. Since $S\ltimes X$ is not contained in the pro-norm of $G$, all subgroups of $T\simeq G/(S\ltimes X)$ are pronormal. Since $\pi(\gamma_3(X))\cap\pi(G/\gamma_3(X))=\pi(X)\cap\pi(G/X)=\emptyset$, we may actually choose $S$ and $T$ such that $G=W\ltimes \gamma_3(X)$, where $W=(T\ltimes S)\ltimes M$ and $X=M\ltimes \gamma_3(X)$. 


The fact that $[M,W]\leq C_M(\gamma_3(X))$ and $W''\leq\operatorname{PAut}_W(\gamma_3(X))$ follows at once from~Lem\-ma~\ref{semidiretto}. In particular, if $\pi([M,W])\cap\pi(C_M(W))=\emptyset$, then we may re-arrange the Sylow subgroups of~$G$ as follows: $G=\big((T\times C_M(W))\ltimes S\big)\ltimes\big([M,W]\times \gamma_3(X)\big)$, where $[M,W]\times\gamma_3(X)$ is Dedekind. This proves the first part of the statement.

\smallskip

Assume $G$ is metabelian. Then $S\ltimes X$ has an ascending $G$-invariant series with elementary abelian factors. 
By Lemma 2.3 of \cite{Goren}, each of these factors is centralized by $T'$ because \hbox{$\pi(T)\cap\pi(SX)=\emptyset$} and $G'$ is abelian. 
Consequently, $T'$ is subnormal in $G$, so it centralizes $SX$. In particular, $T'$ is normal in $G$. Since $T'$ is contained in $T$, so $T'\leq N_{\operatorname{pro}}(G)$ by Lemma \ref{pnormldedek}. Now, $T=U\ltimes \gamma_3(T)$, where $\pi(U)\cap\pi(\gamma_3(T))=\emptyset$ and $U$ is Dedekind, so replacing $T$ by $U$, and $X$ by $X\gamma_3(T)$, we may assume that $T$ is Dedekind. Finally, the Sylow $q$-subgroups of $G$ are~De\-de\-kind for every prime $q$, and hence $G'\leq N_{\operatorname{pro}}(G)$ by Lemma \ref{pnormldedek}.

Now, we prove that $[T,S]=\{1\}$. If $S$ is not Dedekind, then $[S,T]=\{1\}$ by Theorem \ref{semi=diretto}. If $S$ is Dedekind non-abelian, then $T$ centralizes a copy of $Q_8$ contained in $S$ by Lemma \ref{q8lem}. This implies that we can always write  $S=C_S(T)\times [S,T]$. If $x\in C_S(T)$ and $1\neq y\in [S,T]$, then~$\langle xy\rangle$ is not pronormal in~$G$, so $xy\in N_{\operatorname{pro}}(G)$. Since we already know that $[S,T]\leq G'\leq N_{\operatorname{pro}}(G)$, the arbitrariness of $x$ and $y$ yields that $S\leq N_{\operatorname{pro}}(G)$, a contradiction. Thus,~\hbox{$[T,S]=\{1\}$.}



\medskip

\noindent {\it Case \textnormal{($\star$)}}: $[S,M]=\{1\}$.\\Then $S\gamma_3(X)\trianglelefteq G$ and $G/S\gamma_3(X)\simeq T\ltimes M$ has only pronormal subgroups. Re-arranging the~Sy\-low subgroups of $T\ltimes X$ (using also Lemma \ref{q8lem}), we can write $$TX=\big(T\times C_M(T)\big)\ltimes \big([M,T]\times \gamma_3(X)\big),$$ which is a group with only pronormal subgroups. Note that $[S,MT]=\{1\}$, we are done. 

Finally, if $S_1$ is a pronormal subgroup of $S$ that is not contained in the pro-norm of $S$, then it is not contained in the pro-norm of $S\ltimes\gamma_3(X)$, and Lemma \ref{semidiretto} yields that every $S_1$-invariant subgroup of $\gamma_3(X)$ is also $N_S(S_1)$-invariant.

\medskip

Assume conversely that $G$ has the structure described in \textnormal{($\star$)}. Let $H$ be a non-pronormal subgroup of $S\times R$. Then $H$ acts as power automorphisms on $\gamma_3(X)$. We claim that $H$ is contained in the pro-norm of $G$. Let $K$ be any other subgroup of $G$. In order to show that $H$ pronormalizes $K$, we may assume $K\cap [M,T]=\{1\}$, and so even $[M,T]=\{1\}$ by Co\-rol\-la\-ry~3.10 of~\cite{ultimopaper}. Write $K=K_1\ltimes K_2$, where $K_2=\gamma_3(X)\cap X$. Also, we may reduce to $\langle H,K\rangle R\gamma_3(G)$, so now~$K_2$ is normal, and $K_1=K$ is pronormalized by $H$ because it is pronormalized modulo~$\gamma_3(G)$ (see~\cite{ultimopaper},~Co\-rol\-la\-ry~3.10).  Therefore the conditions of Lemma \ref{semidiretto} are satisfied, and we are done.~\end{proof}

\medskip

In our final result, $\overline Z(G)$ denotes the hypercentre of the group $G$.

\begin{theo}\label{q8ipercentro}
Let $G$ be a finite soluble group whose subgroups are pronormal or contained in the pro-norm, and $H$ a subgroup of $G$ which is isomorphic to $Q_8$. If $H\leq\overline{Z}(G)$, then $G=O_2\big(\overline{Z}(G)\big)\times D$ for some $D\leq G$.
\end{theo}
\begin{proof}
By Theorem \ref{pgroups}, we may assume that $G$ is not a $2$-group. Since $G$ is soluble, there are a Sylow $2'$-subgroup $Q$ of $G$ and a Sylow $2$-sub\-\hbox{group~$P$} of $G$ (containing $H$) such that $G=PQ$. Now, $P$ is a quasihamiltonian $2$-group involving $Q_8$, so~\hbox{$P=H\times R$,} where $R$ is an elementary abelian $2$-subgroup. Moreover, since $H\leq\overline{Z}(G)$, we have that $[H,Q]=\{1\}$. Thus, $H$ is normalized by both $P$ and $Q$, and so by $G=PQ$.

Now, let $K$ be any subgroup of $O_2\big(\overline{Z}(G)\big)$ containing $H$. We claim that $K$ is a direct factor of~$G$. First, note that $K$ is normal in $G$ by the previous paragraph because it is generated by copies of $Q_8$. Let $N$ be a minimal normal subgroup of $G$, so $N$ is an elementary abelian~\hbox{$p$-group} for some prime $p$. If $p>2$, then there is $X/N$ such that $G/N=KN/N\times X/N$. Thus, $G=K\times X$ because $K$ is normal in $G$. Suppose that every minimal normal subgroup of $G$ is a $2$-group, put $L=O_2(G)$, and let $S$ be the socle of $L$. Moreover, let $M/L$ be a minimal normal subgroup of~$G/L$, which necessary is a~\hbox{$2'$-sub}\-group. If every subgroup of $S$ is normal in $G$, then $S\leq Z(G)$, so $L\leq \overline Z(G)$ and hence $M=L\times U$ for some normal $2'$-subgroup of $G$, a contradiction. Thus, not every subgroup of $S$ is normal in $G$, and so not every subgroup of $S$ is normalized by $Q$. Thus, there is a minimal normal~\hbox{$2$-sub}\-group~$W$ of order at least $4$ such that $W\cap K=\{1\}$ (by~Ma\-schke's theorem). Finally, by induction, $G/W=KW/W\times D/W$ for some subgroup~$D/W$ of~$G/W$, so $G=K\times D$ and we are done. 
\end{proof}

\section{The infinite soluble case}

The aim of this section is to deal with {\it infinite} (generalized) soluble groups whose subgroups are pronormal or contained in the pro-norm. Most of the results here are essential to prove our main theorem about the general insoluble case (see Theorem \ref{locallycgraded}). We start by showing that some of the results in the previous section in the finite case can be extended to periodic, locally graded groups.

\begin{theo}\label{PLG->LF}
Let $G$ be a locally graded group whose subgroups are pronormal or contained in the pro-norm. If~$G$ is periodic, then $G$ is locally finite.
\end{theo}
\begin{proof}
Clearly, we may assume that $G$ is finitely generated, so we need to prove that $G$ is finite. Let $R$ be finite residual of $G$. If every finite quotient of $G$ is soluble, then every finite quotient of $G$ is soluble of derived length at most $3$ by Corollary \ref{sollength}, so $G/R$ is soluble (of derived length at most $3$) and hence even finite. But then $R=\{1\}$ and we are done. 

It is therefore safe to assume that $G$ contains an insoluble finite quotient $G/N$. Now, $N$ is finitely generated as well, and no finite quotient of $N$ can be insoluble, otherwise every subgroup of $G/N$ would be pronormal and $G/N$ would be soluble, a contradiction. By the previous paragraph, $N$ is finite and so also $G$ is finite. The statement is proved. 
\end{proof}

\begin{rem}
{\rm Theorem \ref{PLG->LF} is improved in Corollary \ref{coramel+}.}
\end{rem}

%
%


We are now going to show that soluble groups whose subgroups are pronormal or contained in the pro-norm satisfy the {subnormal intersection property} (see Theorem \ref{sip}), that is, the intersection of any family of subnormal subgroups is subnormal. In turn, this will provide many structural consequences (see Corollary \ref{cortorsionfree}).

\begin{theo}\label{thserial}
Let $G$ be a group whose subgroups are serial, and pronormal or contained in the pro-norm. Then every subgroup of $G$ is normal or contained in the norm of $G$. In particular, $G$ is quasihamiltonian.
\end{theo}
\begin{proof}
It follows from Theorem 3.1 of \cite{ultimopaper} that for every subgroup $H$ of $G$, the set $P_G(H)$ of all elements of $G$ pronormalizing $H$ coincides with the normalizer $N_G(H)$. Thus, if $H$ is a pronormal subgroup of $G$, then $H$ is normal in $G$. On the other hand, if $H$ is contained in the pro-norm, then $H$ normalizes every subgroup of $G$, so $H$ is actually contained in the norm of~$G$. Thus, $HK=KH$ for all subgroups $H,K$ of~$G$.
\end{proof}

\begin{cor}\label{corLN}
Let $G$ be a locally nilpotent group whose subgroups are pronormal or contained in the pro-norm.  Then every subgroup of $G$ is  normal or contained in the norm of $G$. In particular,~$G$ is quasihamiltonian.
\end{cor}
\begin{proof}
This follows from Theorem \ref{thserial} and the well-known fact that every subgroup of a locally nilpotent group is serial.
\end{proof}

%

\begin{theo}\label{sip}
Let $G$ be a locally graded group whose subgroups are pronormal or contained in the pro-norm. If $H$ is a subnormal subgroup of $G$, then $H$ has defect at most $7$.
\end{theo}
\begin{proof}
Clearly, if $H$ is pronormal in $G$, then $H$ is normal in $G$. Assume $H$ is contained in the pro-norm of $G$. Then all subgroups of $H$ are pronormal in $G$ (and so in $H$), so $H$ is metabelian being locally graded.

Now, $H'$ is an abelian subnormal subgroup of $G$, and so $L=(H')^G$ is locally nilpotent. By~The\-o\-rem~\ref{thserial}, every subgroup of $L$ is normal or contained in the norm of $L$; in particular, every subgroup of $L$ is subnormal of defect at most $2$. 

Let $S$ be a normal subgroup of $L$ such that $H'$ is normal in $S$. Since $L$ and $H'$ are~\hbox{$H$-in}\-va\-riant, we may further require that $S$ is $H$-invariant. We claim that $H$ is subnormal~$HS$ of defect at most $2$. In order to prove this, we may assume $H'=\{1\}$, so $H$ is abelian and subnormal in~$HS$, and by the above $H$ is subnormal of defect at most $2$ in~$HS$. In a similar fashion, we show that~$HS$ is subnormal of defect at most $2$ in $HL$, so $H$ is subnormal of defect at most $4$ in $HL$. Finally, a further application of the same argument shows that $HL$ is subnormal in $(HL)^G$ of defect at most $2$, so $H$ is subnormal in $G$ of defect at most $7$.
\end{proof}

\begin{quest}
What is the smallest possible defect of subnormality in the statement of {\rm Theorem \ref{sip}}?
\end{quest}

\begin{cor}\label{cortorsionfree}
Let $G$ be a soluble group whose subgroups are pronormal or contained in the pro-norm. Then\textnormal: 
\begin{itemize}
    \item[\textnormal{(1)}] If $G$ is finitely generated, then $G$ is finite-by-abelian. Moreover, in this case, $G$ splits over the torsion subgroup.
    \item[\textnormal{(2)}] If $G$ is nilpotent and non-periodic, then $G$ is abelian.
    \item[\textnormal{(3)}] If $G$ is torsion-free, then $G$ is abelian.
    \item[\textnormal{(4)}] If $G$ is non-periodic, then $G$ is metabelian. More precisely, the torsion-subgroup $T$ of $G$ has an abelian normal subgroup $H$ such that $\pi(T)\cap\pi(H)=\emptyset$ and $G/H$ is abelian.
    \item[\textnormal{(5)}] Every non-periodic element $x$ of $G$ has a power $x^n$ centralizing the torsion-subgroup $T$ of $G$, provided that $T$ is abelian. 
\end{itemize}
\end{cor}
\begin{proof}
Clearly, we may always assume that $G$ is finitely generated. Then a combination of~The\-o\-rem~\ref{sip} and \cite{Ro72}, Theorem 10.54, yields that $G$ is finite-by-nilpotent, and so even finite-by-(torsion-free nipotent), being $G$ finitely generated.

Assume $G$ is nilpotent and non-periodic. Then $G$ is quasihamiltonian by Corollary \ref{corLN}, so it satisfies the conditions in \cite{schmidt},~The\-o\-rem~2.4.11. If $G$ is non-abelian, then $G$ has a non-abelian section of the form $\mathcal C_{p^n}\ltimes\mathcal C_{p^m}$, where $m$ is fixed and~$n$ can be chosen arbitrarily large. This contradicts Theorem \ref{pgroups}, and proves (2).

Going back to the general case, we see that $G$ is finite-by-abelian, which proves the first half of (1), and consequently even abelian in case $G$ is torsion-free, thus proving (3).

\smallskip

In order to prove the second half of (1), we may assume that the torsion subgroup $T$ of $G$ is a $p$-group for some prime $p$ (see Theorem \ref{sylowtowerth}). Since $G$ is (central torsion-free)-by-finite, every finite quotient of~$G$ is a quotient of a finite quotient of $G$ in which the centre has $p$-elements of arbitrarily large orders. Thus, Theorem \ref{pgroups} yields that the $p$-subgroups of every quotient of~$G$ are abelian. This implies that there exists a torsion-free central subgroup $Z$ of $G$ such that~\hbox{$G/Z=U/Z\ltimes TZ/Z$} for some $p'$-subgroup $U/Z$. Clearly, $U$ is torsion-free and $G=U\ltimes T$, as desired.

\medskip

\noindent(4)\quad In order to prove the statement, we may assume $G$ is finitely generated, so $G$ is finite-by-(abelian torsion free) by (1). Since $G$ has a torsion-free abelian central subgroup of finite index, we have that $G$ is residually finite-by-(infinite cyclic), so we may actually assume that~$G$ is finite-by-(infinite cyclic). Write $G=\langle w\rangle\ltimes H$, where~$H$ is finite, and $o(w)=\infty$, and set $C=C_{\langle w\rangle}(H)$. It easily follows from Theorem \ref{pgroups} (and Example \ref{exnon} (5)) that the~Sy\-low~\hbox{$p$-sub}\-groups of every quotient of $G$ are abelian (essentially because we can find a larger quotient containing a central~\hbox{$p$-ele}\-ment of large order).

Suppose first that all subgroups of $H$ are pronormal. If $H$ is abelian, then $G$ is metabelian, and we are done. Assume $H$ is not abelian, and write $H=K\ltimes\gamma_3(H)$, where as usual we have $\pi(K)\cap\pi(\gamma_3(H))=\emptyset$, and both $K$ and $\gamma_3(H)$ are abelian. Since the Sylow subgroups of the quotients are abelian, we can write \hbox{$K=C_{\langle w\rangle}(K)\times [\langle w\rangle,K]$} modulo $\gamma_3(H)$. Again, modulo~$\gamma_3(H)$, we have that the non-trivial Sylow subgroups of $[\langle w\rangle,K]\times C$ are not pronormalized by~$w$, so~$K$ is contained in the pro-norm of $G$. On the other hand, $K$ does not pronormalize a cyclic subgroup of $\gamma_3(H)\times\langle w\rangle$, which hence a contradiction. This implies that $[\langle w\rangle,K]=\{1\}$, so $G$ is metabelian.

Now, assume that $H$ does not have only pronormal subgroups. Then $H$ has the structure outlined in Theorem \ref{pronormform}, so we may write $H=(T\ltimes S)\ltimes X$, where $T$, $S$, and $X$ satisfy the conditions in that statement.

We start by deriving some additional information on $H$. Set $L=H\times C$. First, we claim that~$T$, $S$ and $X$ are abelian. By (2), all Sylow subgroups of $G$ are abelian, so $S$ is abelian. Moreover, $T\times\langle a\rangle$ has only pronormal subgroups, so again $T$ is abelian. Let $n$ be the product of all primes in $\pi(H)$, and replace $G$ by $G/C^n$, so $C$ is now finite of order $n$. If $\gamma_3(X)\neq\{1\}$, then $U=\gamma_3(X)\times O_{\pi(\gamma_3(X))}(\langle a\rangle)$ is not contained in the pro-norm, so all subgroups of $G/U$ are pronormal; but $G/U$ has central $p$-subgroups for every prime $p$ in $\pi(G/U)$, so $G/U$ is abelian. If $w$ does not centralize a cyclic subgroup of $H/U$, then $w$ does not normalize a cyclic subgroup of $H/U\times CU/U$, so $S\leq N_{\operatorname{pro}}(G)$, a contradiction.


Now, if $[T,S]\neq\{1\}$, then $T$ does not pronormalize a cyclic subgroup of $S\times C$, so $S$ is contained in the pro-norm of $G$, which is a contradiction. Thus, $[T,S]=\{1\}$. If $w$ does not centralize a cyclic subgroup of $T\times S$, then $w$ does not normalize a cyclic subgroup of $T\times S\times C$, so~\hbox{$T\times S\leq N_{\operatorname{pro}}(G)$,} a contradiction. Consequently, $w$ centralizes $T\times S$, and $G$ is metabelian.

\medskip

\noindent(5)\quad Let $u$ be an element of $T$ that is not centralized by $x$. By (1), there is a power of $x$, say $y=x^n$, centralizing $u$. But then $\langle uy\rangle$ is not pronormal in $\langle uy,x\rangle$, so $uy$ is contained in the pro-norm of $G$, and must consequently act as power automorphisms on $T\times\langle y\rangle$, and so it must centralize $T\times\langle y\rangle$. Thus, $[y,T]=\{1\}$, as desired.
\end{proof}

\begin{rem}
{\rm Example \ref{exyes} (6) shows that we cannot extend Corollary \ref{cortorsionfree} (3) to non-periodic groups, while Example \ref{exyes} (7) shows that non-periodic groups whose subgroups are pronormal or contained in the pro-norm may not even be supersoluble.}
\end{rem}

\begin{cor}\label{boundforsoluble}
Let $G$ be a soluble group whose subgroups are pronormal or contained in the pro-norm. Then $G$ is soluble of length at most $3$.
\end{cor}
\begin{proof}
Clearly, we may assume $G$ is finitely generated. By Corollary \ref{cortorsionfree}, $G$ is polycyclic, and every finite quotient of $G$ is soluble of derived length at most $3$ (see Corollary \ref{sollength}). Since~$G$ is residually finite, $G$ has derived length at most $3$.
\end{proof}

\medskip

Recall that a {\it $T$-group} is a group in which normality is a transitive relation, and that a {\it $\overline T$-group} is a group in which every subgroup is a $T$-group. Soluble $T$-groups have been studied by~Gasch\"utz \cite{Gaschutz} (in the finite case) and Robinson \cite{rob} (in the infinite case). It is a well-known open problem to decide whether a torsion-free locally graded $\overline T$-group is abelian or not. The following result relate the class of torsion-free soluble $\overline T$-groups with the class we are studying.

\begin{cor}
Let $G$ be a torsion-free, locally graded group in which every subgroup is pronormal or contained in the pro-norm. Then $G$ is a $\overline T$-group.
\end{cor}
\begin{proof}
Since the property is inherited by subgroups, we only need to show that $G$ is a~\hbox{$T$-group}. Let $H$ be a subnormal subgroup of $G$. If $H$ is pronormal, then $H$ is normal in $G$. Assume $H$ is contained in the pro-norm of $G$, so $H$ is abelian (being locally graded). Then $H^G$ is locally nilpotent, and so abelian by Corollary \ref{cortorsionfree} (2). Consequently, $H^G\langle x\rangle$ is abelian for every $x\in G$ (by (3) of~Co\-rol\-la\-ry \ref{cortorsionfree}), so $H^G\leq Z(G)$ and $H$ is again normal in $G$.
\end{proof}

    \medskip

Our final result shows that Corollary \ref{cortorsionfree} (5) cannot be much improved.

\begin{ex}
There exists a metabelian non-periodic group with trivial centre whose subgroups are pronormal or contained in the pro-norm.
\end{ex}
\begin{proof}
Let $\{p_i\}_{i\in \mathbb N}$ and $\{q_i\}_{i\in\mathbb N}$ be families of prime numbers such that $q_i$ divides $p_i-1$, and 
$|\{p_i,q_i,p_j,q_j\}|=4$ for $i\neq j$. For each $i\in\mathbb N$, let $\langle a_i\rangle$ and $\langle b_i\rangle$ be groups of order $p_i$ and~$q_i$, respectively. For each $i\in\mathbb N$, let $H_i=\langle b_i\rangle\ltimes\langle a_i\rangle$ be the natural semidirect product. Let~$\langle x\rangle$ be an infinite cyclic group, and define $G=\langle x\rangle\ltimes \operatorname{Dr}_{i\in I} H_i$ by setting $b_i^x=b_ia_i$. Put $H=\operatorname{Dr}_{i\in I} H_i$. Let $X$ be any subgroup of $G$.

Assume $X\not\leq \langle x\rangle H'$. We claim that $X$ is pronormal in $G$. Without loss of generality, we may assume $X\cap H'=\{1\}$, and then even $X\cap H=\{1\}$, so actually $X$ is infinite cyclic. This obviously allows us to assume that $X \le G = \langle x \rangle N$, where $N = \operatorname{Dr}_{i=1}^nH_i$ for $n$ sufficiently large, and $\pi(\langle x\rangle)=\pi(N)$. Now, $Z(G) \le \langle x \rangle$, and $G/Z(G)$ is finite, so we may assume that $Z(G) = \{1\}$, and hence that $G$ is finite. 
Setting $z= xa_1^{-1} \ldots a_n^{-1}$, we have that $G = \langle z \rangle \times N$.  By Lemma \ref{semidiretto}, every subgroup of $G$ is pronormal or contained in the pro-norm. Moreover, by~Lem\-ma~\ref{pnormldedek}, $N_{\operatorname{pro}}(G) = \langle z \rangle \times N'=\langle x\rangle N'$ (note that the pro-norm is maximal in $G$ and $G$ has not only pronormal subgroups), so $X$ is pronormal in $G$.

Now, suppose $X\leq \langle x\rangle H'$. Let $W$ be a subgroup of $G$. We need to show that $X$ pronormalizes~$W$. By the previous paragraph, we may assume that $W\leq \langle x\rangle H'$. But $\langle x\rangle H'$ is abelian, and so we are done.
\end{proof}

\section{The insoluble case}

The aim of this section is to give a complete description of finite insoluble groups whose subgroups are pronormal or contained in the pro-norm. This is the content of Theorem \ref{locallycgraded}, which is actually proved in a large class of insoluble (arbitrary) groups, although the main bulk of the proof is certainly the finite case. The layout of the proof is as follows: first we are going to characterize the finite simple groups (see Theorem \ref{classificationsimple}), then the direct products of finite simple groups with arbitrary finite groups, and finally we are going exploit this information to prove the main result. 

We start by proving that certain finite simple groups satisfy our property. Most of the information about the structure of the finite simple groups involved here comes from  \cite{atlas} and~\cite{King}.

\begin{theo}\label{j1}
All subgroups of the first Janko group $J_1$ are pronormal or contained in the pro-norm.
\end{theo}
\begin{proof}
Put $G=J_1$. Let $H$ be a non-pronormal subgroup of $G$. We need to show that $H$ is contained in the pro-norm of $G$. By \cite{perReferee4}, every non-nilpotent subgroup of $G$ is pronormal. The nilpotent subgroups of $G$ are cyclic or elementary abelian $2$-groups (see \cite{atlas}). Cyclic subgroups of odd order are pronormal by \cite{perReferee4}, Corollary 2.4, and \cite{vdovin}. Those of orders $6$ and $10$ are also pronormal: they contain a pronormal Sylow $3$- or $5$-subgroup $P$ of $G$, respectively, and $N_G(P)\simeq\operatorname{Sym}(3)\times\operatorname{Dih}(10)$ has all its subgroups pronormal (see \cite{KS} and \cite{perReferee4}, Lemma 2.1). Since the Sylow $2$-subgroups have order $8$ and are pronormal, the only possible non-pronormal subgroups have order $2$ or $4$, and so without loss of generality we may assume that $H=\langle x\rangle$ is cyclic of order $2$. 

Now, let $K$ be any subgroup of $G$, and put $L=\langle K,K^{x}\rangle$. The proof will be complete if we can prove that $K$ is conjugate to $K^{x}$ in $L$. To this aim, we may obviously assume that $L$ is a proper subgroup of $G$. Furthermore, $x$ normalizes $L$ because $o(x)=2$, so $\langle x\rangle L$ is in fact a proper subgroup of $G$ (by the simplicity of $G$), and hence it is contained in a maximal subgroup $M$ of $G$. Then $M$ is isomorphic to one of the following groups:
\begin{itemize}
    \item[(1)] $\operatorname{PSL}(2,11)$;\quad (2) $(\mathcal C_3\ltimes\mathcal C_7)\ltimes\mathcal C_2^3$;\quad (3) $\mathcal C_2\times\operatorname{Alt}(5)$;\quad (4) $\mathcal C_{6}\ltimes \mathcal C_{19}$;
    \item[(5)] $\mathcal C_{10}\ltimes\mathcal C_{11}$;\quad (6) $\operatorname{Sym}(3)\times\operatorname{Dih}(10)$;\quad (7) $\mathcal C_6\ltimes\mathcal C_7$.
\end{itemize}

\noindent The groups of types (4)--(7) have all their subgroups pronormal (see \cite{KS}). In type (2), every involution belongs to the normal elementary abelian $2$-subgroup of $M$, and hence to the pro-norm of $M$ by Lemma \ref{pnormldedek}. In type (3), every involution belongs to the pro-norm of $M$, as a direct check of the subgroups of $\mathcal C_2\times\operatorname{Alt}(5)$ shows (see also \cite{GAP}). Finally, type (1) is $\operatorname{PSL}(2,11)$, whose involutions belong to its pro-norm by Theorem \ref{finalpsl}. Thus $x$ pronormalizes $K$ in every case, and the statement is proved.
\end{proof}

\medskip

\begin{theo}\label{psl2}
Let $p$ be a prime and set $G=\operatorname{PSL}(2,2^p)$. Then every subgroup of $G$ is pronormal or contained in the pro-norm.
\end{theo}
\begin{proof}
By \cite{perReferee4}, every non-nilpotent subgroup of $G$ is pronormal. Its nilpotent subgroups are elementary abelian $2$-groups or cyclic groups of odd order (see \cite{King}). Each such cyclic group $U$ is pronormal: choose a subgroup $P\leq U$ of prime order, which is pronormal by \cite{perReferee4}, Corollary 2.4. Since $N_G(P)$ is dihedral of order $2(2^p-1)$ or $2(2^p+1)$, all its subgroups are pronormal, and the conclusion follows from \cite{perReferee4}, Lemma 2.1. Sylow $2$-subgroups are pronormal as well. Thus the only non-pronormal subgroups are elementary abelian $2$-groups. Let $H=\langle x\rangle$ be any cyclic subgroup of order $2$, and let $K$ be any subgroup of $G$. As in the proof of Theorem \ref{j1}, we may assume that $H$ normalizes $\langle K,K^x\rangle$, so~$\langle H,K\rangle$ is a proper subgroup of $G$. By \cite{King}, each maximal subgroup $M$ is either of the form $E\ltimes C$, with $E$ a normal elementary abelian $2$-group and $C$ cyclic of odd order, or is dihedral of order $2(2^p-1)$ or $2(2^p+1)$. In the first case, every involution lies in $E\leq N_{\operatorname{pro}}(M)$ by Lemma \ref{pnormldedek}; in the second, all subgroups of $M$ are pronormal (see \cite{KS}). Thus $x$ pronormalizes $K$ and proves the result. 
\end{proof}

\begin{lem}\label{chat}
Let $q$ be an odd prime power with $q \equiv_8 5$,
\[
y_x=\begin{pmatrix}1&x\\[2pt]0&1\end{pmatrix}\in \operatorname{SL}(2,q),\; x\neq0,
\]
and $t\in \operatorname{PSL}(2,q)$ an involution. Then 
$\langle t,\,t^{y_x}\rangle$ is non-abelian.
\end{lem}
\begin{proof}
Since $t$ is an involution in $\operatorname{PSL}(2,q)$ and $q\equiv 5\pmod 8$, every lift 
$A\in \operatorname{SL}(2,q)$ of $t$ has trace $\operatorname{Tr}(A)=0$; hence we may write
\[
A=\begin{pmatrix}a & b\\[2pt] c & -a\end{pmatrix}, \qquad  
a^2+bc=-1,
\]
where the last identity comes from $\det(A)=1$. Clearly,
\[
A^{y_x}=y_x^{-1}Ay_x=
\begin{pmatrix}
a-cx & b+2ax-cx^2\\[2pt]
c   & -a+cx
\end{pmatrix}.
\]

Denote $P=AA^{y_x}$ and $Q=A^{y_x}A$. A direct matrix multiplication gives
\[
P=\begin{pmatrix}
a^2-acx+bc & 2a^2x-acx^2+bcx\\[2pt]
-c^2x      & a^2+acx+bc-c^2x^2
\end{pmatrix},
\]
\[
Q=\begin{pmatrix}
a^2+acx+bc-c^2x^2 & -2a^2x+acx^2-bcx\\[2pt]
c^2x              & a^2-acx+bc
\end{pmatrix}.
\]

We now analyse when the classes of $A$ and $A^{y_x}$ commute in $\operatorname{PSL}(2,q)$.
Since the centre of $\operatorname{SL}(2,q)$ is $\{\pm I_2\}$, commuting in $\operatorname{PSL}(2,q)$ means that
\[
AA^{y_x}=\varepsilon\, A^{y_x}A,
\qquad 
\varepsilon\in\{\pm 1\}.
\]

\noindent{\it Case $\varepsilon=+1$}\quad
The equation $P=Q$ implies, by comparing the $(2,1)$-entries, that
\[
-c^2x=c^2x,
\]
hence $2c^2x=0$. Since $q$ is odd and $x\neq 0$, we obtain $c=0$.
Substituting $c=0$ into the $(1,2)$-entries we get
\[
2a^2x=-2a^2x,
\]
and therefore $4a^2x=0$. Again, since $q$ is odd and $x\neq 0$, this yields
$a=0$. But then $a^2+bc=0$, contradicting the identity $a^2+bc=-1$.
Thus no such $A$ exists, and $A$ and $A^{y_x}$ do not commute in
$\operatorname{SL}(2,q)$.

\medskip

\noindent{\it Case $\varepsilon=-1$}\quad
The equation $P=-Q$ is equivalent to $P+Q=0$, and we compute
\[
P+Q=\begin{pmatrix}
2a^2+2bc-c^2x^2 & 0\\[2pt]
0 & 2a^2+2bc-c^2x^2
\end{pmatrix}
= (2a^2+2bc-c^2x^2)\,I_2.
\]
Hence $P+Q=0$ if and only if
\[
2a^2+2bc-c^2x^2=0.
\]
Using $a^2+bc=-1$, we substitute $a^2=-1-bc$ and obtain
\[
2(-1-bc)+2bc-c^2x^2=0
\quad\Longleftrightarrow\quad
-2-c^2x^2=0
\quad\Longleftrightarrow\quad
c^2x^2=-2.
\]

Thus $t$ and $t^{y_x}$ commute in $\operatorname{PSL}(2,q)$ if and only if 
$c^2=-2/x^2$. Since $x^2$ is a square, this holds if and only if $-2$ is a square
in $\mathbb F_q$, which is false for $q\equiv_8 5$.
Therefore $\langle t,t^{y_x}\rangle$ is non-abelian.
\end{proof}

\begin{theo}\label{5mproof}
Let $G = \operatorname{PSL}(2,q)$, where $q=5^m$, and either $m=1$ or $m$ is an odd prime. Let~$H$ be a subgroup of~$G$. Then $H$ is pronormal in $G$ provided that one of the following conditions holds\textnormal:
    \begin{itemize}
        \item[\textnormal{(1)}] $H$ is cyclic with $|H|$ dividing $\frac{q\pm1}{2}$ and $|H|\neq 2$.
        \item[\textnormal{(2)}] $H$ is dihedral with $|H|$ dividing $q\pm1$.
        \item[\textnormal{(3)}] $H\simeq\operatorname{Alt}(4)$ or $H\simeq\operatorname{Alt}(5)$.
        \item[\textnormal{(4)}] $H$ is a Hall subgroup of $G$.
        \item[\textnormal{(5)}] $S\leq H\leq N_G(S)$, where $S$ is a Sylow $5$-subgroup of $G$, without any restriction on $|H/S|$.
    \end{itemize}

\noindent Moreover, $H$ is contained in the pro-norm of $G$ provided that\textnormal:
    \begin{itemize}
        \item[\textnormal{(6)}] $H$ is cyclic of order $2$ or of order $5$.
    \end{itemize}

\smallskip

\noindent In particular, every subgroup of $G$ is pronormal or contained in the pro-norm.
\end{theo}
\begin{proof}
(1)\quad 
Since $q$ is odd, $(q+1)/2$ and $(q-1)/2$ are coprime, and hence the Sylow~\hbox{$\pi$-sub}\-groups of $G$ are cyclic whenever $\pi$ is a set of odd primes dividing $q+1$ or $q-1$. Then a well-known result of John Rose (see for example \cite{perReferee4}, Corollary 2.4) shows that $H$ is pronormal in $G$ whenever $|H|$ is odd.

Now, assume that $|H| = 2m$ where $m$ is odd, and let $L$ be the subgroup of $H$ of order $m$. By the previous paragraph, $L$ is pronormal in $G$. Moreover, $C_G(L)$ is cyclic (the centralizer of any non-trivial element of odd order prime to $5$ in $\operatorname{PSL}(2,q)$ is cyclic; see \cite{King}). Let $g\in G$. Then there  exists $x \in \langle L,L^g \rangle  \leq  \langle H,H^g \rangle$ such that $L^g = L^x$, and consequently $H^x,H^g \leq C_G(L^x) = C_G(L^g)$. But $C_G(L^g)$ is cyclic, and $|H^x|=|H^g|$, so  $H^g=H^x$. Thus, $H$ is pronormal in $G$.

    \smallskip

\noindent(2)\quad Write $H=C\ltimes L$, where $|C|=2$ and $L$ is  cyclic. Of course, we may assume that $4\neq|L|>2$. Set $R=O_{2'}(L)$, and let $g\in G$. We need to show that~$H$ and $H^g$ are conjugate in $\langle H,H^g\rangle$. If $\langle H,H^g\rangle=G$, then we are good. Assume not, and let $x\in\langle H,H^g\rangle$ be such that $R^x=R^g$. Clearly, it is enough to show that~$H^x$ and $H^g$ are conjugate in $\langle H^x,H^g\rangle\leq \langle H,H^g\rangle$, and so that~$H$ and~$H^{gx^{-1}}$ are conjugate in $\langle H,H^{gx^{-1}}\rangle$. We may therefore replace $g$ by $gx^{-1}$ obtaining that $R=R^g$. Then $\langle H,g\rangle\leq N_G(R)<G$, so $\langle H,g\rangle$ is dihedral. Since $q\equiv_8 5$, the cyclic subgroup of index $2$ in $N_G(R)$ has order not divisible by $4$. Thus all subgroups of $N_G(R)$ are pronormal (see \cite{KS}), and $H$ and $H^g$ are conjugate in $\langle H,H^g\rangle$, as required.

\medskip

\noindent(3)\quad If $H\simeq\operatorname{Alt}(5)$, then $H$ is maximal. Suppose $H\simeq\operatorname{Alt}(4)$ and let $g\in G$. If~$\langle H,H^g\rangle=G$, then we are done. Thus, $\langle H,H^g\rangle$ is either isomorphic to $\operatorname{Alt}(4)$ or $\operatorname{Alt}(5)$. But in both cases, $H$ is conjugate to $H^g$ inside the generated subgroup.

\medskip

\noindent(4)\quad For Sylow subgroups this is obvious, while for Hall subgroups in general follows from~\cite{vdovin}.

\medskip

\noindent(5)\quad Since $N_G(S)/S$ is cyclic, every subgroup $H$ with $S\leq H\leq N_G(S)$ is normal in $N_G(S)$. As $S$ is pronormal in $G$, \cite{perReferee4}, Lemma 2.1, implies that $H$ is pronormal in $G$.

\medskip

\noindent(6)\quad Let $K$ be any subgroup of $G$, and write $H=\langle x\rangle$. If $o(x)=2$, then $H$ normalizes $\langle K,K^x\rangle$, so we may assume that $\langle H,K\rangle$ is a proper subgroup of $G$. But every maximal subgroup of $G$ is easily seen to have the involutions contained in the pro-norm.

Assume $o(x)=5$. By (1)--(5) and Lemma \ref{chat}, we may assume that $K\leq N_G(S)$ for some~Sy\-low~\hbox{$5$-sub}\-group $S$ of $G$, and $|K|=5^\ell\cdot 2^\varepsilon$, where $\varepsilon\in\{0,1\}$. Let~\hbox{$R=O_5(K)$.} If $[R,R^x]=\{1\}$, then $R^x\leq Z(\langle R,R^x\rangle)\cap Z(\langle R^x,R^{x^2}\rangle)$, and so $R^x\leq S$. Repeating this argument yields that $K^H$ is contained in $S$ and so $R$ is normalized by $H$. Then $H\leq N_G(S)$. But $H$ is contained in the pro-norm of $N_G(S)$, and so we are done.

Suppose $[R,R^x]\neq\{1\}$ and $\langle K,K^x\rangle<G$. Then $\langle K,K^x\rangle\simeq\operatorname{Alt}(5)$, and either $|K|=5$ or $|K|=10$. In both cases, there is only one conjugacy class of subgroups of that order, so we are done again.
\end{proof}

\begin{theo}\label{finalpsl}
Let $G=\operatorname{PSL}(2,q)$, where $q>5$ is a prime with $q\equiv_83,5$. Then every subgroup of $G$ is pronormal or contained in the pro-norm.
\end{theo}
\begin{proof}
The proof is essentially the same as that of Theorem \ref{5mproof} but easier, so it is omitted.~\end{proof}

\begin{theo}\label{classificationsimple}
Let $G$ be a non-abelian finite simple group. Then every subgroups of $G$ is pronormal or contained in the pro-norm if and only if $G$ is isomorphic to one of the following types of groups\textnormal:
\begin{itemize}
    \item[\textnormal{(1)}] $\operatorname{PSL}(2,5^m)$, where either $m=1$ or $m$ is an odd prime\textnormal; 
    \item[\textnormal{(2)}] $\operatorname{PSL}(2,q)$, where $q$ is a prime with $q>5$ and $q\equiv_8 3,5$\textnormal;
    \item[\textnormal{(3)}] $\operatorname{PSL}(2,2^m)$, where $m$ is a prime\textnormal;
    \item[\textnormal{(4)}] the first Janko group $J_1$.
\end{itemize}
\end{theo}
\begin{proof}
The sufficiency of the conditions is obtained by combining Theorems \ref{j1}, \ref{psl2}, \ref{5mproof}, and \ref{finalpsl}. Let us turn to the necessity. Let $S$ be a Sylow $2$-subgroup of $G$. Then $S$ has nilpotency class at most $2$, by~The\-o\-rem \ref{pgroups}. Assume $S$ has class $2$. Then $G$ is isomorphic to one of the following groups (see \cite{Goren2}): $\operatorname{PSL}(2,q)$, with $q\equiv_{16} 7,9$, $\operatorname{Alt}(7)$, $\operatorname{Sz}(2^n)$, $\operatorname{PSL}(3,2^n)$, $\operatorname{PSp}(4,2^n)$, $U(3,2^n)$, where $n\geq2$. However, in any of the first three cases, the Sylow $2$-subgroups have a section that is isomorphic with $\operatorname{Dih}(8)$, which contradicts Theorem \ref{pgroups}.

Suppose  $G =\operatorname{PSp(4,2^n) = Sp(4,2^n)}$ (see \cite{Huppert}, p.231). Now, $G$ contains a~\hbox{$2$-sub}\-group~$Q$ such that $Q' =  Z(Q) = \operatorname{Fratt}(Q)$ is elementary abelian of order at least $2^n$ (see \cite{Wilson}, p.63), contradicting Theorem \ref{pgroups}.

Suppose $G=U(3,2^n)$. Let $\widetilde G=\operatorname{SU}(3,2^n)$ and $Z=Z(\widetilde G)$, so that $G\simeq\widetilde G/Z$ (see \cite{Wilson}, p.66). The group $Z$ has odd order $d=(2^n+1,3)$, and hence the Sylow $2$-subgroups of $G$ and $\widetilde G$ are isomorphic. It follows that $G$ contains a $2$-subgroup $Q$ of order $2^{3n}$ such that $Z(Q) = Q' = \operatorname{Fratt}(Q)$ (see \cite{Wilson}, pg. 67), contradicting Theorem \ref{pgroups}. 

Thus, $S$ is abelian and $G$ is isomorphic to one of the following groups (see for example \cite{Goren}, p.485): $\operatorname{PSL}(2,q)$, where $q >3$, $q \equiv_8 3,5$ or $q = 2^n$, $J_1$, $^2G_2(3^{2n+1})$. But the~Sy\-low~\hbox{$3$-sub}\-groups of $^2G_2(3^{2n+1})$ are of nilpotency class $3$ (see for example \cite{Busch}, p.5), so they can be ruled out.


Write $q=q_0^m$, and let $Q$ be a Sylow $q_0$-subgroup of $G$ (this is an elementary abelian $q$-subgroup). Then $N_G(Q)=C\ltimes Q$ for some cyclic subgroup $C$ irreducibly acting on $Q$. If $m$ is divided by at least $2$ primes, then there are proper subgroups $D_1\leq D_2$ of $C$ and proper subgroups $R_1\leq R_2$ of $Q$ such that $R_2^{D_2}=R_2$, $R_1^{D_1}=R_1$ and $R_1^{D_2}\neq R_1$. Then $D_2R_2$ is not pronormalized by any generator of $C$, and, at the same time, it is not even contained in the pro-norm of $G$ because $D_1$ is not. Thus, either $m=1$ or $m$ is a prime.

Suppose $m$ is prime and $q$ is odd. Since $q\equiv_83,5$ and $q$ is odd, so $m>2$. If $q_0=3$, then $G$ contains a cyclic subgroup of order $3$ as a non-pronormal subgroup. But these subgroups are conjugate and $G$ also contains a copy of $\operatorname{Alt}(4)$, so it contains a subgroup of order $3$ that is not pronormal nor contained in the pro-norm of $G$. Thus $q_0\neq3$.

Assume $q\equiv_83$ and $m>1$. Since $m$ is odd, we necessarily have $q_0\equiv_83$ (just check all the possibilities for $q_0$). This implies that $-2=a^2$ for some $a\in \operatorname{GF}(q_0)$. Let $y$ and $x$ be respectively the images of the matrices $$
\begin{pmatrix}
1 & a\\
0 & 1
\end{pmatrix}\qquad\textnormal{and}\qquad \begin{pmatrix}
0 & 1\\
-1 & 0
\end{pmatrix}
$$ in $G$. It is easy to see that $\langle x,x^y\rangle$ has order $4$, so $y$ does not pronormalize $\langle x\rangle$. On the other hand, $\langle y\rangle$ is certainly not pronormal because $m>1$. This rules out the case $q\equiv_83$ and \hbox{$m>1$.}

Finally, assume $p>5$, $q\equiv_85$ and $m>1$. We shall choose an involution $x$ and an element $y$ of $G$ represented by the matrices $$
X=\begin{pmatrix}
a & b\\
c & -a
\end{pmatrix}
$$ and $$
Y=\begin{pmatrix}
t & 0\\
0 & t^{-1}
\end{pmatrix}\qquad (t\in\mathbb F_p)
$$ in $\operatorname{SL}(2,q)$, where $a^2+bc=-1$, $t\neq0$, and $t^4\neq1$. The values of $a,b,c,t$ will be specified below. We wish to find a necessary and sufficient condition for $\langle x,x^y\rangle$ to be of order $4$. First, we compute $$
X^Y=\begin{pmatrix}
a & bt^{-2}\\
ct^2 & -a
\end{pmatrix}
$$ and then $$XX^Y=\begin{pmatrix}
a^2+bct^2 & ab(t^{-2}-1)\\
ac(1-t^2) & a^2+bct^{-2}
\end{pmatrix},\qquad X^YX=\begin{pmatrix}
a^2+bct^{-2} & ab(1-t^{-2})\\
ac(t^2-1) & a^2+bct^2
\end{pmatrix}.$$ Requiring that $XX^Y=X^YX$ is equivalent to requiring that $$\begin{cases}
ab(t^{-2}-1)=ab(1-t^{-2})\\
ac(1-t^2)=ac(t^2-1)\\
a^2+bct^2=a^2+bct^{-2}
\end{cases}$$ But $t^4\neq 1$, so the previous system is equivalent to $ab=ac=bc=0$, which gives $b=c=0$ in virtue of $a^2+bc=-1$, and so that $y\in C_G(x)$, a contradiction. Requiring that $XX^Y=-X^YX$ is equivalent to requiring that $2a^2+bc(t^2+t^{-2})=0$. Replacing $bc=-1-a^2$ in the previous equation gives $a^2(2-(t^2+t^{-2}))=t^2+t^{-2}$. Clearly, 
$2-t^2-t^{-2}=0$ if and only if $(t^2-1)=0$, so $$a^2=\frac{t^2+t^{-2}}{2-(t^2+t^{-2})}.$$ The previous equation can be solved for $a$ if and only if $-(t^4+1)$ is a square in $\mathbb F_q$. Since $m$ is odd, this is equivalent to requiring that $-(t^4+1)$ is a square in $\mathbb F_p$.

Let $\chi$ denote the quadratic character of $\mathbb{F}_p$, extended by
$\chi(0)=0$.
Thus $\chi(a)=1$ if $a$ is a non-zero square, $\chi(a)=-1$ if $a$ is a non-square. Consider the character sum
\[
S := \sum_{t \in \mathbb{F}_p} \chi\bigl(-(t^4+1)\bigr).
\]
Using multiplicativity of $\chi$, we have
\[
S = \chi(-1)\sum_{t \in \mathbb{F}_p} \chi(t^4+1).
\]

The polynomial $f(t)=t^4+1$ has degree $4$ and distinct roots.
Therefore the Weil bound for multiplicative character sums applies (see for example \cite{weil}):
\[
\left| \sum_{t \in \mathbb{F}_p} \chi(t^4+1) \right|
\le (4-1)\sqrt{p} = 3\sqrt{p}.
\]
Hence
\[
|S| \le 3\sqrt{p}.
\]

We argue by contradiction.
Assume that for every $t\in\mathbb F_p$ the element $-(t^4+1)$ is not a square. Since $t^4+1=0$ has never a root in $\mathbb F_p$ because $-1$ is not a square, then $|S|=p$ and consequently $p<3\sqrt{p}$, which gives $p\leq 5$, a contradiction. Thus, there is $t\in\mathbb F_p$ such that $-(t^4+1)$ is a square in $\mathbb F_p$. Clearly, $t\neq 0$ and $t^2\neq1$ because $-2$ is not a square in $\mathbb F_p$ because $p\equiv_85$. Choose $a$ satisfying $a^2=-(t^4+1)/(t^2-1)^2$, and set $b=1$ and $c=-1-a^2$. Then $X$ represents an involution $x$, $y\notin C_G(x)$, and the calculation above gives $|\langle x,x^y\rangle|=4$.

\medskip

Now, let $\ell$ be the order of $y$. By our assumption on $t$, $y$ is not contained in the pro-norm of~$G$. On the other hand, since $m>1$, there is a subgroup of order $p\ell$ containing $y$ and which is not pronormal in the normalizer of a Sylow $p$-subgroup of $G$. This contradicts the hypotheses and completes the proof of the statement.
\end{proof}



\begin{cor}\label{noinfinitesimple}
There is no infinite locally graded periodic simple group whose subgroups are pronormal or contained in the pro-norm.
\end{cor}
\begin{proof}
Let $G$ be an infinite locally graded periodic simple group whose subgroups are pronormal or contained in the pro-norm. By Theorem \ref{PLG->LF}, $G$ is locally finite.  By~The\-o\-rem~\ref{classificationsimple}, we may find a Kegel cover $\mathcal K$ of $G$ in such a way that every Kegel factor is linear of bounded degree (see \cite{KW}). Since the finite subgroups of $G$ are isomorphic with subgroups of some Kegel factor of $\mathcal K$, so $G$ is linear by a well-known theorem of Mal'cev. But then $G$ is the union of its finite simple subgroups (see \cite{KW}, Proposition~4.6), a fact which is clearly impossible because of the classification provided in~The\-o\-rem~\ref{classificationsimple}.
\end{proof}

\medskip

\begin{cor}\label{corpgl}
Let $G=\operatorname{PGL}(n,p^m)$, where $p$ is a prime, $n\geq2$, $m\geq1$. Every subgroup of $G$ is pronormal or contained in the pro-norm if and only if $p=2$, and either $m$ is a prime, or $m=1$ and $n=2$.
\end{cor}
\begin{proof}
This is immediate from Theorem \ref{classificationsimple} and \cite{King}, Corollary 2.3.
\end{proof}

 \medskip

 The following easy result allows us to extend Theorem \ref{classificationsimple} to the corresponding special linear groups.

\begin{lem}\label{centralextension}
Let $G$ be a periodic group, $Z$ a central subgroup of $G$, $U$ a subgroup of $G$, and $x\in P_{G/Z}(UZ/Z)$. Then $x\in P_G(U)$ provided that $\pi(U)\cap\pi(Z)=\emptyset$, or $Z\leq U$.
\end{lem}
\begin{proof}
The statement is clear if $Z\leq U$. Assume $\pi(U)\cap\pi(Z)=\emptyset$. Then there exists an element $y\in \langle U, U^x \rangle Z$
such that $U^y\times Z=U^yZ = U^xZ=U^x\times Z$. Thus, $U^y=U^x$. But $y=hz^m$ for some $h\in\langle U,U^x\rangle$ and $z\in Z$, so $U^h=U^x$. 
\end{proof}

\begin{cor}\label{corcentralextension}
Let $G$ be a periodic group, and $Z$ a central subgroup of $G$ of order $2$ such that every subgroup of $G/Z$ is pronormal or contained in the pro-norm of $G/Z$. If the Sylow $2$-subgroups of~$G$ are generalized quaternion,  then every subgroup of $G$ is  pronormal or contained in the pro-norm of $G$. \textnormal(In particular, the Sylow $2$-subgroups must be quaternion groups of order $8$.\textnormal)
\end{cor}

\begin{theo}\label{autpsl22p}
Let $G=\operatorname{Aut}(\operatorname{PSL}(2,2^p))$, where $p$ is a prime. Every subgroup of $G$ is pronormal or contained in the pro-norm if and only if $p>3$. In this case, the subgroups of order $2$ are contained in the pro-norm.
\end{theo}
\begin{proof}
If $p=2$, then $G\simeq\operatorname{Sym}(5)$, which contains a cyclic subgroup of order $2$ that is neither pronormal nor contained in the pro-norm. If $p=3$, then we can easily find a subgroup of the form $\mathcal C_2\ltimes (\mathcal C_3\times\mathcal C_3)$ where the action of $\mathcal C_2$ is not the same of the two $\mathcal C_3$ factors (this is the centralizer of the Frobenius automorphism). Since every subgroup of order $2$ of $\operatorname{PSL}(2,2^p)$ is not pronormal, we have again found a cyclic subgroup of order $2$ that is neither pronormal nor contained in the pro-norm. 

Assume $p>3$. Write $G=K\ltimes\operatorname{PSL}(2,2^p)$, where $|K|=p$, and note that $(p,|\operatorname{PSL}(2,2^p)|)=1$. All pronormal subgroups of $\operatorname{PSL}(2,2^p)$ are again pronormal in $G$ because their isomorphism classes in $\operatorname{PSL}(2,2^p)$ coincides with their conjugacy classes in $\operatorname{PSL}(2,2^p)$. If $H$ is a subgroup of $G$ whose order is divided by~$p$, then we can assume $H=K\ltimes R$ for some subgroup $R$ of~$\operatorname{PSL}(2,2^p)$. If $R$ is pronormal in $\operatorname{PSL}(2,2^p)$, then~$H$ is pronormal in $G$, being the product of two pronormal subgroups, one of which normalizes the other. Thus, we may assume that $R$ is an elementary abelian~\hbox{$2$-group} of order $2^m$, where $0<m<p$.

Let $g\in G$. In order to show that $H$ is pronormal in $G$, we may assume that $H<L=\langle H,H^g\rangle<G$. Write $L=K\ltimes S$ for some $R<S<\operatorname{PSL}(2,2^p)$. Of course, we may also assume that $H$ is not a Hall subgroup of $L$, so $L$ is contained in the normalizer of a~Sy\-low~\hbox{$2$-sub}\-group~$S$ of $G$. Clearly, we may further assume $g\in\operatorname{PSL}(2,2^p)$.

Now, $\langle R,R^g\rangle\cap\langle R^g,R^{g^2}\rangle\neq\{1\}$, so $R$, $R^g$ and $R^{g^2}$ all lie in $S$ (see for example \cite{Huppert},~The\-o\-rem~8.2, p.200). Repeating this argument, we see that $g$ normalizes $S$. On the other hand, there exists $u\in\langle K, K^g \rangle \le L$ such that $K^u = K^g$, and so $gu^{-1} \in N_G(K)=K\times \operatorname{SL}(2,2)\simeq\mathcal C_p\times\operatorname{Sym}(3)$. Since $g,u$ are elements of $\operatorname{PSL}(2,2^p)$ normalizing $S$, then $gu^{-1}$ has order $2$ or $3$, but it cannot be $3$ because $p$ is odd and $3$ does not divide $2^p-1$. Thus, $gu^{-1}$ has order $2$ and it must necessarily be contained in $S$, which implies that $$(K \ltimes R)^g = K^g \ltimes R^g = K^u \ltimes R^g = K^u \ltimes R^{(gu^{-1})u} = K^u \ltimes R^u = (K \ltimes R)^u.$$ Therefore $H$ is pronormal in $G$.


Finally, let $x$ be any element of order $2$ and $H$ any subgroup of $G$. In order to show that $x$ is contained in the pro-norm of $G$, we may assume $H$ is not pronormal, so $H$ is an elementary abelian $2$-group. But $x$ normalizes $\langle H,H^x\rangle$, so $\langle x,H\rangle$ is an abelian $2$-subgroup and consequently $H^x=H$.
\end{proof}

\begin{theo}\label{aut5^p}
Let $G$ be the subgroup of index $2$ of $\operatorname{Aut}(\operatorname{PSL}(2,5^p))$, where $p$ is an odd prime. Then~$G$ has a subgroup that is neither pronormal nor contained in the pro-norm.
\end{theo}
\begin{proof}
Write $G = K \ltimes U$, where $U=\operatorname{PSL}(2,5^p)$, and $|K| = p$. Note that $$C_G(K)=K\times \operatorname{PSL}(2,5) \simeq \mathcal C_p\times \operatorname{Alt}(5).$$
    Let $A$ be a subgroup of $\operatorname{PSL}(2,5)$ that is isomorphic to $\operatorname{Alt}(4)$, and $C$ a subgroup of order $2$ of~$A$. 
    The subgroup $H = K \times C$ is not pronormal in $K \times A$, so it is not pronormal in $G$. 
    On the other hand, $K$ is not contained in the pro-norm of $G$, because it is not contained in the pro-norm of~$K\ltimes V$, where $V$ is the image in $U$ of the group of all unitriangular matrices of~$\operatorname{SL}(2,5^p)$.
\end{proof}

\begin{theo}\label{J1direct}
Let $G=H\times K$ be a locally graded group, where $K\simeq J_1$. Then all subgroups of $G$ are pronormal or contained in the pro-norm if and only if the following conditions are satisfied\textnormal:
\begin{itemize}
    \item[\textnormal{(1)}] $H$ is periodic, and all its subgroups are pronormal.
    \item[\textnormal{(2)}] $\pi(H)\cap\pi(K)\subseteq\{2\}$.
\end{itemize}
\end{theo}
\begin{proof}
First, we deal with the necessity. Clearly, $K$ has subgroups that are not contained in the pro-norm, so every subgroup of $H\simeq G/K$ must be pronormal. Suppose that there is an odd prime $p$ in $\pi(H)\cap\pi(K)$. Then $G$ has a subgroup of the form $C\ltimes (\mathcal C_p\times\mathcal C_p)$, where $C$ is a subgroup of order $2$ of $K$ acting in distinct ways on the factors (see the maximal subgroups of~$J_1$ in the proof of Theorem \ref{j1}). It follows that $C$ is not contained in the pro-norm of $G$. But~$C$ is not even pronormal in $K$, and so we have a contradiction. This proves (2). Finally, (2) yields that $H$ does not contain any infinite cyclic subgroup, and so that $H$ is periodic, proving~(1).

\medskip

Now, we deal with the sufficiency. Write $H=R\ltimes S$, where $S$ is abelian with all subgroups normal in $H$, $R$ is Dedekind, $\pi(R)\cap\pi(S)=\emptyset$, $2\not\in\pi(S)$. Let $U$ be a subgroup of $G$. We need to show that $U$ is pronormal in $G$ or is contained in the pro-norm of $G$. 

\medskip

\noindent{\it Step \textnormal(a\textnormal)}\quad If $U=U_H\times U_K$, where $U_H\leq H$ and $U_K\leq K$ with $|U_K|\neq2,4$, then $U$ is pronormal in~$G$. 
(The proof of this is obvious because $U_H$ and $U_K$ are pronormal in $H$ and $K$, respectively.)


\medskip

\noindent {\it Step \textnormal(b\textnormal)}\quad Suppose that $|UH/H|\neq 2,4$ (so $UH/H$ is pronormal in $G/H$), and that~$U$ cannot be written as $U_H\times U_K$, where $U_H\leq H$ and $U_K\leq K$. We claim that $U$ is pronormal in~$G$.\\
\indent Let $g\in G$. Since \hbox{$\pi(S)\cap\pi(K)=\emptyset$,} Corollary 3.16 of \cite{ultimopaper} shows that we may assume that $S=\{1\}$ and $H$ is a Dedekind $2$-group. Consequently, we may further assume that $U\cap H=\{1\}$, $H$ is an elementary abelian \hbox{$2$-group}, and $g\in K$. In particular, $U$ is finite, so, replacing $G$ by~\hbox{$\langle U,K,g\rangle$,} we may as well assume that $G$ is finite. We need to divide the proof according to the possibilities for~$U$. Suppose first that $19$ divides $|U|$, so $U\cap K$ contains a normal subgroup~$W$ of order $19$, and we may assume $U$ is contained in a subgroup of the form $(\langle a\rangle\times B)\ltimes W$, where $|B|=3$, and $a=xy$, where $x\in H$ and $y\in K$ are involutions. If $W^g\neq W$, then $\langle U,U^g\rangle\geq \langle W,W^g\rangle =K\ni g$. Thus, we may assume $W^g=W$, so $g\in\langle W,B,y\rangle$, and hence we easily see that $U$ and $U^g$ are both contained in $U_1=\langle W,B,a\rangle$. But $U$ and $U^g$ are Hall subgroups of $U_1$ and so they are conjugate in the subgroup they generate. Thus, we may assume that $19$ does not divide $|U|$.

Suppose that $U$ is insoluble. Then $U\cap K=K$, $U\cap K\simeq\operatorname{PSL}(2,11)$, $U\cap K\simeq\operatorname{Alt}(5)$, or $U\cap K\simeq\mathcal C_2\times \operatorname{Alt}(5)$. In both the first and second cases, we obtain the contradiction that $U=(U\cap H)\times (U\cap K)$ (because the elements of $G$ act on $U\cap K$ as elements of $K$). In the fourth case, we may assume that $g\in N_K(U\cap K)=U\cap K$, and so we are done. Assume that~\hbox{$U\cap K\simeq\operatorname{Alt}(5)$,} so $U\simeq \langle xy\rangle\times\operatorname{Alt}(5)$, where $x\in H$ and $y\in K$ are involutions. If $g$ normalizes $U\cap K$, then $g\in\langle y,U\cap K\rangle$, so $g$ centralizes $xy$ and we are done. Thus, $g$ does not normalize $U\cap K$. If $\langle U\cap K,(U\cap K)^g\rangle\simeq\operatorname{PSL}(2,11)$, then $g\in K\leq \langle xy, U\cap K,(U\cap K)^g\rangle$, so we are done. Finally, if $\langle U\cap K,(U\cap K)^g\rangle=\langle y,U\cap K\rangle$, then obviously $U\cap K=(U\cap K)^g$, a contradiction. Thus, we may assume $U$ is soluble.

Now, the case in which $11$ divides $|U|$ can be dealt as that in which $19$ divides~$|U|$. Note only that $U$ (containing a dihedral subgroup of order $22$) is not isomorphic to any subgroup of $\operatorname{PSL}(2,11)$, so (using the notation of the previous case) if $\langle W,W^g\rangle\simeq\operatorname{PSL}(2,11)$, then $K\leq\langle xy, W,W^g\rangle\leq \langle U,U^g\rangle$. It is therefore possible to assume that $11$ does not divide the order of~$|U|$.

Next, suppose $5$ divides the order of $|U|$. The previous assumptions imply that $U$ has a normal subgroup $W$ of order $5$. First note that there exists $z\in\langle U,U^g\rangle\cap K$ such that $W^z=W^g$. In order to prove that $U$ is conjugate to $U^g$ in $\langle U,U^g\rangle$, we only need to prove that $U^z$ and~$U^g$ are conjugate in $\langle U^z,U^g\rangle\leq\langle U,U^g\rangle$, and so that $U$ and $U^{gz^{-1}}$ are conjugate in $\langle U,U^{gz^{-1}}\rangle$. Replacing $g$ with $gz^{-1}$, we may assume $W=W^g$.  Now, $N_K(W)=\big(\langle k_1\rangle\ltimes \langle v\rangle\big)\times\big(\langle k_2\rangle\ltimes W\big),$ where $o(k_1)=o(k_2)=2$, $o(v)=3$, $v^{k_1}=v^{-1}$, and $w^{k_2}=w^{-1}$ for $w\in W$. Consequently,  $$U\leq\big(\langle h_1k_1\rangle\ltimes \langle v\rangle\big)\times\big(\langle h_2k_2\rangle\ltimes W\big),$$ where $h_1,h_2\in H$. Since $g\in N_K(W)$, so $\langle U,U^g\rangle\leq\langle h_1,h_2\rangle\times N_K(W)$. But all subgroups of~\hbox{$\langle h_1,h_2\rangle\times N_K(W)$} are pronormal in it, and so we are done. Thus, we may assume that $5$ does not divide the order of $U$. 

A similar argument applies when $U$ has a normal subgroup of order $7$ or $3$, so we may assume $U$ has no such subgroups. (Note that this argument could have also been employed in dealing with the cases in which $19$ (or $11$) divides the order of $|U|$ but we chose to exhibit a different argument.)

Now, assume that $7$ divides $|U|$. By the previous assumptions, we have that $U$ contains a normal Sylow $2$-subgroup $P$ of order $8$ (which is elementary abelian) and an element $x$ of order~$7$ acting irreducibly on $P$. But then $P$ must lies inside $K$, which in turns implies that $U\leq K$, a contradiction. Therefore $7$ does not divide the order of~$|U|$.

Finally, assume that $3$ divides $|U|$. Then either $U\simeq\operatorname{Alt}(4)$, in which case we obtain the contradiction $U\leq K$, or $U=LM$, where $M\simeq\operatorname{Alt}(4)$, $L$ is cyclic of order $6$, and $L$ normalizes~$M$. But we already know that $L$ and $M$ are pronormal subgroups of $G$, so $U$ is pronormal as well.

\medskip

\noindent {\it Step \textnormal(c\textnormal)}\quad If $g$ is any involution of $K$, then $H\times \langle g\rangle$ is contained in the pro-norm of $G$.\\
\indent Let $U$ be a subgroup of $G$, and let $x\in H\times\langle g\rangle$. We need to show that $U$ and $U^x$ are conjugate in $\langle U,U^x\rangle$. To this aim we may assume $U$ is not pronormal, so by Steps (a) and (b), $UH/H$ is a (finite) $2$-subgroup of~$G/H$. Again, Corollary 3.16 of \cite{ultimopaper} shows that we may assume $H$ to be a (Dedekind) $2$-group, so $x$ is a $2$-element. Consequently, we may further assume $U\cap H=\{1\}$, so~$U$ is a finite elementary abelian~\hbox{$2$-group.} Now, replacing $G$ by $\langle U,x,K\rangle$, we may further assume that $G$ is finite. Write $x=yz$, where $y\in H$ and $z\in\langle g\rangle\leq K$. Then the action of $x$ on~$U$ coincides with that of $z$, so we may assume that $x\in K$, that $H$ is an elementary abelian~\hbox{$2$-group}, and that $U$ is not contained in $K$.

Clearly, we may further assume that $\langle U,U^x\rangle H/H<G/H$, so also $\langle U,x\rangle H/H<G/H$ because~$xH$ has order $2$ (and so normalizes $\langle U,U^x\rangle$). It follows that $\langle U,x\rangle\leq H\times M$, where $M$ is a maximal subgroup of $K$ with no normal Sylow $2$-subgroup (otherwise $\langle U,x\rangle$ is abelian). But in all these cases, we easily see that the involutions are contained in the pro-norm. The claim is proved.

\medskip

By combining Steps (a), (b), and (c), the statement is proved.
\end{proof}

\begin{theo}\label{H5mod8}
Let $G=H\times K$ be a locally graded group, where $K\simeq \operatorname{PSL}(2,q)$ for some prime $q\geq 5$ with $q \equiv_{8} 5$. Then all subgroups of $G$ are pronormal or contained in the pro-norm if and only if the following conditions are satisfied\textnormal:
\begin{itemize}
    \item[\textnormal{(1)}] $H$ is periodic, and all its subgroups are pronormal.
    \item[\textnormal{(2)}] $\pi(H)\cap\pi(K)\subseteq\{2\}$.
\end{itemize}
\end{theo}
\begin{proof}
The necessity is proved as in Theorem \ref{J1direct} using the list in Theorem 2.1 of \cite{King}.

\medskip

Now, we deal with the sufficiency. Write $H=R\ltimes S$, where $S$ is abelian with all subgroups normal in $H$, $R$ is Dedekind, $\pi(R)\cap\pi(S)=\emptyset$, $2\not\in\pi(S)$. Let $U$ be a subgroup of $G$. We need to show that $U$ is pronormal in $G$ or is contained in the pro-norm of $G$. 

\medskip

\noindent{\it Step \textnormal(a\textnormal)}\quad If $U=U_H\times U_K$, where $U_H\leq H$ and $U_K\leq K$ with $|U_K|\neq2,$ then $U$ is pronormal in~$G$. (The proof of this is obvious because $U_H$ and $U_K$ are pronormal in $H$ and $K$, respectively.)

\medskip

\noindent {\it Step \textnormal(b\textnormal)}\quad Suppose that $|UH/H|$ is not a power of $2$ (so $UH/H$ is pronormal in $G/H$), and that~$U$ cannot be written as $U_H\times U_K$, where $U_H\leq H$ and $U_K\leq K$. We claim that $U$ is pronormal in~$G$.\\
\indent Let $g\in G$. Since \hbox{$\pi(S)\cap\pi(K)=\emptyset$,} Corollary 3.16 of \cite{ultimopaper} shows that we may assume that $S=\{1\}$ and $H$ is a Dedekind $2$-group. Consequently, we may further assume that $U\cap H=\{1\}$, $H$ is an elementary abelian \hbox{$2$-group}, and $g\in K$. In particular, $U$ is finite, so, replacing $G$ by~\hbox{$\langle U,K,g\rangle$,} we may as well assume that $G$ is finite. We need to divide the proof according to the possibilities for~$U$. 

If $U$ is insoluble, then $U\simeq\operatorname{Alt}(5)$ must be contained in $K$, contradicting our assumptions. Thus, $U$ is solvable. Let $N$ be a minimal normal subgroup of $U$. If $|N|=4$, then $U\simeq\operatorname{Alt}(4)$, and so $N\leq K$ as well as $U\leq K$, a contradiction. 

Thus, we may assume that $|N|\neq2$, so $N$ is contained in a cyclic Sylow subgroup of $K$. Consequently, there is~\hbox{$z\in \langle U,U^g\rangle\cap K$} such that $N^z=N^g$. In order to prove that $U$ is conjugate to $U^g$ in $\langle U,U^g\rangle$, we only need to prove that $U^z$ and $U^g$ are conjugate in $\langle U^z,U^g\rangle\leq\langle U,U^g\rangle$, and so that $U$ and~$U^{gz^{-1}}$ are conjugate in $\langle U,U^{gz^{-1}}\rangle$. Replacing $g$ with $gz^{-1}$, we may assume $N=N^g$. 

Now, we can write $U=\langle x_1y_1\rangle\ltimes\big(\langle x_2y_2\rangle\times R\big)$, where $x_1,x_2\in H$ and $y_1,y_2\in K$ have order at most $2$, and $R$ is cyclic of odd order. Then $\langle U,g\rangle\leq \langle x_1\rangle\times\langle x_2\rangle\times N_K(N)$, where $N_K(N)=\langle y_1,y_2,S\rangle$ with $S$ cyclic of odd order containing $R$. Clearly, all subgroups of $L$ are pronormal, so $U$ is conjugate to $U^g$ in $\langle U,U^g\rangle$, and the claim is proved.

\medskip

\noindent {\it Step \textnormal(c\textnormal)}\quad If $g$ is any involution of $K$, then $H\times \langle g\rangle$ is contained in the pro-norm of $G$.\\
\indent (This is proved exactly as the corresponding step in the proof of Theorem \ref{J1direct}.)

\medskip

By combining Steps (a), (b), and (c), the result follows. 
\end{proof}

\begin{cor}\label{corsl2,q}
Let $G=H\times K$ be a locally graded group, where $K\simeq \operatorname{SL}(2,q)$ for some prime $q\geq 5$ with $q \equiv_{8} 5$. Then all subgroups of $G$ are pronormal or contained in the pro-norm if and only if the following conditions are satisfied\textnormal:
\begin{itemize}
    \item[\textnormal{(1)}] $H$ is periodic, and all its subgroups are pronormal.
    \item[\textnormal{(2)}] $\pi(H)\cap\pi(K)\subseteq\{2\}$.
    \item[\textnormal{(3)}] The Sylow $2$-subgroups of $H$ are elementary abelian.
\end{itemize}
\end{cor}
\begin{proof}
For the sufficiency, it is enough to apply Theorems \ref{H5mod8} and Example \ref{exnon} (5), having noted that $G/Z(K)$ is obviously locally graded. Let's deal with the necessity.  By Theorem \ref{H5mod8}, we know that all subgroups of $G/Z(K)$ are  pronormal or contained in the pro-norm. Let $x\in G$ and $U\leq G$ with $x$ pronormalizing $UZ(K)/Z(K)$. We claim that $x$ pronormalizes $U$: this is enough to complete the proof.  If $2\not\in\pi(U)$ or $U\geq Z(K)$, then we may apply~Lem\-ma~\ref{centralextension}. Assume $2\in\pi(U)$ and $U\cap Z(K)=\{1\}$, so the Sylow $2$-subgroups of $U$ are elementary abelian and are consequently centralized by $K$. By~Co\-rol\-la\-ry~3.16 of \cite{ultimopaper}, we may assume that $H$ is abelian, so the Sylow $2$-subgroups of $U$ are normal in $G$. By factoring them out, we see that $x$ pronormalizes $U$.
\end{proof}

\begin{theo}\label{H5^pmod8}
Let $G=H\times K$ be a locally graded group, where $K\simeq \operatorname{PSL}(2,5^p)$ for some odd prime~$p$. Then all subgroups of $G$ are pronormal or contained in the pro-norm if and only if the following conditions are satisfied\textnormal:
\begin{itemize}
    \item[\textnormal{(1)}] $H$ is periodic, and all its subgroups are pronormal.
    \item[\textnormal{(2)}] $\pi(H)\cap\pi(K)\subseteq\{2\}$.
\end{itemize}
\end{theo}
\begin{proof}
The necessity follows at once from Theorem \ref{H5mod8}. Now, to the sufficiency. Write $H=R\ltimes S$, where $S$ is abelian with all subgroups normal in $H$, $R$ is Dedekind, $\pi(R)\cap\pi(S)=\emptyset$, $2\not\in\pi(S)$. Let $U$ be a subgroup of $G$. We need to show that $U$ is pronormal in $G$ or is contained in the pro-norm of $G$. 

\medskip

\noindent{\it Step \textnormal(a\textnormal)}\quad If $U=U_H\times U_K$, where $U_H\leq H$ and $U_K\leq K$ with $\pi(U_K)\not\subseteq \{2,5\}$, then $U$ is pronormal in~$G$. (The proof of this is obvious because by Theorem \ref{5mproof} $U_H$ and $U_K$ are pronormal in~$H$ and~$K$, respectively.)

\medskip

\noindent {\it Step \textnormal(b\textnormal)}\quad Suppose that $\pi(UH/H)\not\subseteq\{2,5\}$ (so $UH/H$ is pronormal in $G/H$ by Theorem \ref{5mproof}), and that~$U$ cannot be written as $U_H\times U_K$, where $U_H\leq H$ and $U_K\leq K$. We claim that $U$ is pronormal in~$G$.\\
\indent Let $g\in G$. Since \hbox{$\pi(S)\cap\pi(K)=\emptyset$,} Corollary 3.16 of \cite{ultimopaper} shows that we may assume that $S=\{1\}$ and $H$ is a Dedekind $2$-group. Consequently, we may further assume that $U\cap H=\{1\}$, $H$ is an elementary abelian \hbox{$2$-group}, and $g\in K$. In particular, $U$ is finite, so, replacing $G$ by~\hbox{$\langle U,K,g\rangle$,} we may as well assume that $G$ is finite. We need to divide the proof according to the possibilities for~$U$. 

If $U$ is insoluble, then $U$ is simple and so it must be contained in $K$, contradicting our assumptions. Thus, $U$ is solvable. Let $N$ be a minimal normal subgroup of $U$. If $|N|=4$, then~$U\simeq\operatorname{Alt}(4)$, and so $N\leq K$ as well as $U\leq K$, a contradiction. 

Thus, we may assume that $|N|\neq2$, so $N$ is either a Sylow $5$-subgroup of $G$, or it is contained in a cyclic Sylow subgroup of $K$. Consequently, there is~\hbox{$z\in \langle U,U^g\rangle\cap K$} such that $N^z=N^g$. In order to prove that $U$ is conjugate to $U^g$ in $\langle U,U^g\rangle$, we only need to prove that $U^z$ and~$U^g$ are conjugate in $\langle U^z,U^g\rangle\leq\langle U,U^g\rangle$, and so that $U$ and~$U^{gz^{-1}}$ are conjugate in $\langle U,U^{gz^{-1}}\rangle$. Replacing $g$ with $gz^{-1}$, we may assume $N=N^g$. 

Suppose $5$ does not divide $|N|$. In this case we can write $U=\langle x_1y_1\rangle\ltimes\big(\langle x_2y_2\rangle\times R\big)$, where $x_1,x_2\in H$ and $y_1,y_2\in K$ have order at most $2$, and $R$ is cyclic of odd order. Then $\langle U,g\rangle\leq L =\langle x_1\rangle\times\langle x_2\rangle\times N_K(N)$,
where $N_K(N)=\langle y_1,y_2,S\rangle$ with $S$ cyclic of odd order containing~$R$. Clearly, all subgroups of $L$ are pronormal, so $U$ is conjugate to $U^g$ in $\langle U,U^g\rangle$, and we are done.

If $5$ divides $|N|$, then $N$ is a Sylow $5$-subgroup of $G$, and $U=C\ltimes N$, where $C$ is a cyclic subgroup such that $2\mid|C|\mid (5^p-1)$. Now, $N$ is obviously pronormal and $C$ is pronormal by the previous cases, so $U=CN$ is pronormal as well. This completes Step (b).

\medskip

\noindent {\it Step \textnormal(c\textnormal)}\quad If $g$ is any element of $K$ of order $2$ or $5$, then $H\times \langle g\rangle$ is contained in the pro-norm of~$G$.\\
\indent Let $U$ be a subgroup of $G$, and let $x\in H\times\langle g\rangle$. We need to show that $U$ and $U^x$ are conjugate in $\langle U,U^x\rangle$. To this aim we may assume $U$ is not pronormal, so by Steps (a) and~(b), $\pi(UH/H)\subseteq\{2,5\}$. Again,~Co\-rol\-la\-ry~3.16 of \cite{ultimopaper} shows that we may assume $H$ to be a (Dedekind) $2$-group. Consequently, we may further assume $U\cap H=\{1\}$, so $UK/K$ is an elementary abelian~\hbox{$2$-group} and hence $o(x)\in\{2,5,10\}$. Now, $U$ is finite, so, replacing $G$ by~\hbox{$\langle U,K,g\rangle$,} we may as well assume that $G$ is finite. Write $x=yz$, where $y\in H$ and $z\in K$. Then the action of $x$ on $U$ coincides with that of $z$, and so we may assume that $x\in K$ (in particular, $o(x)\in\{2,5\}$), that $H$ is an elementary abelian $2$-group, and that $U$ is not contained in $K$.

Clearly, we may further assume that $\langle U,U^x\rangle H/H<G/H$. If $\langle U,U^x\rangle H/H$ is insoluble, then it is simple, so $\langle U,U^x\rangle\leq K$ and hence applying The\-o\-rem~\ref{5mproof} we are done. Thus, $\langle U,U^x\rangle H/H$ is solvable. 

If $\langle U,x\rangle H/H$ is a proper subgroup of $G/H$, then $\langle U,x\rangle\leq H\times M$, where~$M$ is a maximal subgroup of $K$ with no normal Sylow $2$-subgroup (otherwise $\langle U,x\rangle$ is abelian). In all these cases, we see that the elements of order $2$ and $5$ are contained in the pro-norm. On the other hand, if we assume $\langle U,x\rangle H=G$, then the argument employed in point (6) of the proof of~The\-o\-rem~\ref{5mproof} yields that $\langle U,x\rangle H/H$ is contained in the normalizer of a Sylow $5$-subgroup of $G/H$, giving a contradiction. Step (c) is proved.


\medskip

By combining Steps (a), (b), and (c), the result follows. 
\end{proof}

\begin{cor}\label{corsl5q}
Let $G=H\times K$ be a locally graded group, where $K\simeq \operatorname{SL}(2,5^p)$ for some odd prime~$p$. Then all subgroups of $G$ are pronormal or contained in the pro-norm if and only if the following conditions are satisfied\textnormal:
\begin{itemize}
    \item[\textnormal{(1)}] $H$ is periodic, and all its subgroups are pronormal.
    \item[\textnormal{(2)}] $\pi(H)\cap\pi(K)\subseteq\{2\}$.
    \item[\textnormal{(3)}] The Sylow $2$-subgroups of $H$ are elementary abelian.
\end{itemize}
\end{cor}
\begin{proof}
The proof is essentially the same as that of Corollary \ref{corsl2,q}.
\end{proof}

\begin{lem}\label{chat-q3}
Let $q$ be an odd prime power with $q \equiv_8 3$,
\[
y_x=\begin{pmatrix}1&x\\[2pt]0&1\end{pmatrix}\in \operatorname{SL}(2,q),\; x\neq0.
\] Then there exists an involution $t\in \operatorname{PSL}(2,q)$ such that
$\langle t,\,t^{y_x}\rangle$ is abelian.
\end{lem}
\begin{proof}
Since $q\equiv 3\pmod 8$, every involution $t\in \operatorname{PSL}(2,q)$
lifts to a matrix
\[
A=\begin{pmatrix}a & b\\[2pt] c & -a\end{pmatrix}\in \operatorname{SL}(2,q),
\qquad a^2+bc=-1,
\]
with $\operatorname{Tr}(A)=0$. Set $A^{y_x}=y_x^{-1}Ay_x$. As shown in the proof of Lemma~\ref{chat}, the classes of $A$ and~$A^{y_x}$ commute in $\operatorname{PSL}(2,q)$ if and only if
\[
AA^{y_x}=-A^{y_x}A,
\]
which is equivalent to
\[
c^2x^2=-2.
\]

Since $q\equiv 3\pmod 8$, the element $-2$ is a square in $\mathbb F_q$.
Fix $\lambda\in\mathbb F_q$ such that $\lambda^2=-2$, and define
\[
c=\lambda/x.
\]
Then $c^2x^2=-2$.

Now choose
\[
a=0
\qquad\text{and}\qquad
b=-\frac{1}{c}.
\]
With this choice we have
\[
a^2+bc = 0 + b c = -1,
\]
so $A\in \operatorname{SL}(2,q)$.
Hence $A$ lifts an involution $t\in \operatorname{PSL}(2,q)$ such that
$t$ and $t^{y_x}$ commute in~$\operatorname{PSL}(2,q)$. Therefore $\langle t,t^{y_x}\rangle$ is abelian.
\end{proof}

\begin{theo}
Let $G=H\times K$ be a locally graded group, where $K\simeq \operatorname{PSL}(2,q)$ for some odd prime $q$ with $q \equiv_{8} 3$. Then all subgroups of $G$ are pronormal or contained in the pro-norm if and only if the following conditions are satisfied\textnormal:
\begin{itemize}
    \item[\textnormal{(1)}] $H$ is periodic, and all its subgroups are pronormal.
    \item[\textnormal{(2)}] $\pi(H)\cap\pi(K)\subseteq\{2\}$.
\end{itemize}
\end{theo}
\begin{proof}
First, we deal with the necessity. Clearly, $K$ has subgroups that are not contained in the pro-norm, so every subgroup of $H\simeq G/K$ must be pronormal. Suppose that there is an prime $p \notin \{2,q\}$ in $\pi(H)\cap\pi(K)$. Then $G$ has a subgroup of the form \hbox{$C\ltimes (\mathcal C_p\times\mathcal C_p)$,} where $C$ is a subgroup of order $2$ of $K$ inverting one $\mathcal C_p$ factor and centralizing the other (see~\cite{King}, The\-o\-rem~2.1). It follows that $C$ is not contained in the pro-norm of $G$. But~$C$ is not even pronormal in $K$, and so we have a contradiction. 

Now, we deal with the prime $q$. Let $x\in H$ and $y\in K$ be elements of order $q$, and let~\hbox{$C=\langle xy\rangle$.} Then $C$ is not pronormal because there exist elements of $K$ normalizing but not centralizing~$\langle y\rangle$. On the other hand, Lemma \ref{chat-q3} shows that $y$ does not belong to the pro-norm of $K$, so $xy$ does not belong to the pro-norm of $G$. Thus, $C$ is neither pronormal nor contained in the pro-norm, a contradiction. This proves (2). Finally, (2) yields that $H$ does not contain any infinite cyclic subgroup, and so that $H$ is periodic, proving~(1).
    
\medskip

Now, we deal with the sufficiency. Write $H=R\ltimes S$, where $S$ is abelian with all subgroups normal in $H$, $R$ is Dedekind, $\pi(R)\cap\pi(S)=\emptyset$, $2\not\in\pi(S)$. Let $U$ be a subgroup of $G$.
We need to show that $U$ is pronormal in $G$ or is contained in the pro-norm of $G$. 

\medskip

\noindent{\it Step \textnormal(a\textnormal)}\quad If $U=U_H\times U_K$, where $U_H\leq H$ and $U_K\leq K$ with $|U_K|\neq2,$ then $U$ is pronormal in~$G$. (The proof of this is obvious because by Theorem \ref{5mproof} $U_H$ and $U_K$ are pronormal in~$H$ and~$K$, respectively.)

\medskip

\noindent {\it Step \textnormal(b\textnormal)}\quad Suppose that $\pi(UH/H)\not\subseteq\{2\}$ (so $UH/H$ is pronormal in $G/H$), and that~$U$ cannot be written as $U_H\times U_K$, where $U_H\leq H$ and $U_K\leq K$. We claim that $U$ is pronormal in~$G$.\\
\indent This can be proved as the corresponding point in the proof of Theorem \ref{H5mod8}.

\medskip

\noindent {\it Step \textnormal(c\textnormal)}\quad If $g$ is any element of $K$ of order $2$, then $H\times \langle g\rangle$ is contained in the pro-norm of~$G$.\\
\indent This can be proved as the corresponding point in the proof of Theorem \ref{J1direct}.

\medskip

By combining Steps (a), (b), and (c), the result follows. 
\end{proof}

\begin{cor}\label{corsl2q}
Let $G=H\times K$ be a locally graded group, where $K\simeq \operatorname{SL}(2,q)$ for some odd prime $q$ with $q \equiv_{8} 3$. Then all subgroups of $G$ are pronormal or contained in the pro-norm if and only if the following conditions are satisfied\textnormal:
\begin{itemize}
    \item[\textnormal{(1)}] $H$ is periodic, and all its subgroups are pronormal.
    \item[\textnormal{(2)}] $\pi(H)\cap\pi(K)\subseteq\{2\}$.
    \item[\textnormal{(3)}] The Sylow $2$-subgroups of $H$ are elementary abelian.
\end{itemize}
\end{cor}
\begin{proof}
The proof is the same as that of Corollary \ref{corsl2,q}.
\end{proof}

\begin{theo}\label{plusautpsl22}
Let $G=H\times K$ be a locally graded group, where $K\simeq \operatorname{PSL}(2,2^p)$ for some prime $p\geq3$. Then all subgroups of $G$ are pronormal or contained in the pro-norm if and only if the following conditions are satisfied\textnormal: 
    \begin{itemize}
        \item[\textnormal{(1)}] $H$ is periodic, and all its subgroups are pronormal.
        \item[\textnormal{(2)}] $\pi(H) \cap \pi(K) \subseteq \{2\}$. 
    \end{itemize}
\end{theo}
\begin{proof}
First, we deal with the necessity. Clearly, $K$ has subgroups that are not contained in the pro-norm, so every subgroup of $H\simeq G/K$ must be pronormal. Suppose that there is an odd prime $p$ in $\pi(H)\cap\pi(K)$. Then $G$ has a subgroup of the form $C\ltimes (\mathcal C_p\times\mathcal C_p)$, where $C$ is a subgroup of order $2$ of $K$ inverting one $\mathcal C_p$ factor and centralizing the other (see~\cite{King}, The\-o\-rem~2.1). It follows that $C$ is not contained in the pro-norm of $G$. But $C$ is not even pronormal in $K$ (otherwise being central in the normalizer of a Sylow $2$-subgroup of $K$), and so we have a contradiction. This proves (2).  Finally, (2) yields that $H$ does not contain any infinite cyclic subgroup, and so that $H$ is periodic, proving~(1).

\medskip

Now, we deal with the sufficiency. Write $H=R\ltimes S$, where $S$ is abelian with all subgroups normal in $H$, $R$ is Dedekind, $\pi(R)\cap\pi(S)=\emptyset$, $2\not\in\pi(S)$. Let $U$ be a subgroup of $G$. We need to show that $U$ is pronormal in $G$ or is contained in the pro-norm of $G$. 

\medskip

\noindent{\it Step \textnormal(a\textnormal)}\quad If $U=U_H\times U_K$, where $U_H\leq H$, $U_K\leq K$, and $|U_K|$ is not a power of $2$, then $U$ is pronormal in~$G$. (The proof of this is obvious because $U_H$ and $U_K$ are pronormal in~$H$ and~$K$, respectively.)

\medskip

\noindent {\it Step \textnormal(b\textnormal)}\quad Suppose that $|UH/H|$ is not a power of $2$ (so $UH/H$ is pronormal in $G/H$), and that~$U$ cannot be written as $U_H\times U_K$, where $U_H\leq H$ and $U_K\leq K$. We claim that $U$ is pronormal in~$G$.\\
\indent Let $g\in G$. Since \hbox{$\pi(S)\cap\pi(K)=\emptyset$,} Corollary 3.16 of \cite{ultimopaper} shows that we may assume that $S=\{1\}$ and $H$ is a Dedekind $2$-group. Consequently, we may further assume that $U\cap H=\{1\}$, $H$ is an elementary abelian \hbox{$2$-group}, and $g\in K$. In particular, $U$ is finite, so, replacing $G$ by~\hbox{$\langle U,K,g\rangle$,} we may as well assume that $G$ is finite. We need to divide the proof according to the possibilities for~$U$. 

Since $K$ is a minimal simple group, $U$ is solvable. Let $N$ be a minimal normal subgroup of $U$. If $|N|$ is a power of $2$, then $|N|=2^p$ and $N$ is irreducibly acted upon by an odd-order element, so $N\leq K$ as well as $U\leq K$, a contradiction. 

Thus, we may assume that $|N|$ is not a power of $2$, so $N$ is contained in a cyclic Sylow subgroup of $K$. Consequently, there is~\hbox{$z\in \langle U,U^g\rangle\cap K$} such that $N^z=N^g$. In order to prove that $U$ is conjugate to $U^g$ in~\hbox{$\langle U,U^g\rangle$,} we only need to prove that $U^z$ and $U^g$ are conjugate in~\hbox{$\langle U^z,U^g\rangle\leq\langle U,U^g\rangle$,} and so that $U$ and~$U^{gz^{-1}}$ are conjugate in $\langle U,U^{gz^{-1}}\rangle$. Replacing $g$ with~$gz^{-1}$, we may assume $N=N^g$. 

Now, we can write $U=\langle xy\rangle\ltimes  R$, where $x\in H$ and $y\in K$ have order at most $2$, and $R$ is cyclic of odd order. Then $\langle U,g\rangle\leq L=\langle x\rangle\times N_K(N)$, where $N_K(N)=\langle y,S\rangle$ with $S$ cyclic of odd order containing $R$. Clearly, all subgroups of $L$ are pronormal, so $U$ is conjugate to $U^g$ in~\hbox{$\langle U,U^g\rangle$,} and the claim is proved.

\medskip

\noindent {\it Step \textnormal(c\textnormal)}\quad If $g$ is any involution of $K$, then $H\times \langle g\rangle$ is contained in the pro-norm of $G$.\\
\indent This can be proved as the corresponding point in the proof of Theorem \ref{J1direct}.

\medskip

By combining Steps (a), (b), and (c), the statement is proved.
\end{proof}

\begin{theo}\label{damigliorare}
Let $p>3$ be a prime, and let $G = C \ltimes (H \times K)$, where $C=\langle c\rangle$ is cyclic of~\hbox{$p$-power} order, $K=\operatorname{PSL}(2,2^p)$, $H\neq\{1\}$ is locally graded, and $[C,K]\neq\{1\}$ \textnormal(in particular, $C/C^{p}\!\ltimes\! K\!\simeq\! \operatorname{Aut}(\operatorname{PSL}(2,2^p))$\textnormal). Then all subgroups of $G$ are pronormal or contained in the pro-norm if and only if the following conditions are satisfied\textnormal: 
    \begin{itemize}
        \item[\textnormal{(1)}] $H$ is periodic, and all subgroups of $C \ltimes H$ are pronormal. 
        \item[\textnormal{(2)}] $H\trianglelefteq G$, $K\trianglelefteq G$.
        \item[\textnormal{(3)}] $\pi(H) \cap \pi(K) \subseteq \{2\}$.
    \end{itemize}
\end{theo}
\begin{proof}
Note that $p$ does not divide the order of $K$.  We start dealing with the necessity. It follows from Theorem \ref{plusautpsl22} that $\pi(H)\cap\pi(K)\subseteq\{2\}$, and that $H$ is a periodic group with only pronormal subgroups. In particular, $H$ is soluble, and so both $H$ and $K$ are characteristic in~\hbox{$H\times K$.} Consequently, $H$ and $K$ are normal subgroups of $G$. Since $K$ has subgroups of order~$2$ that are not contained in the pro-norm, $G/K\simeq C\ltimes H$ has only pronormal subgroups. The necessity is proved. 

        \medskip

Now, we deal with the sufficiency. Set $L:=C\ltimes H$ and write $L=R\ltimes S$, where $S$ is abelian with all subgroups normal in $L$, $R$ is Dedekind, $\pi(R)\cap\pi(S)=\emptyset$, $2\not\in\pi(S)$, $S\leq L'$, and $C\leq R$. Let $U$ be a subgroup of $G$.
We need to show that $U$ is pronormal in $G$ or is contained in the pro-norm of $G$. 

\medskip

\noindent{\it Step \textnormal(a\textnormal)}\quad If $U=U_H\times U_K$, where $U_H\leq C^p\ltimes H$, $U_K\leq K$, and $|U_K|$ is not a power of $2$, then $U$ is pronormal in~$G$.
 
 Let $g = c^nhk\in G$, where $n\in\mathbb N_0$, $h \in H$, and $k \in K$. Then $$(U_H\times U_K)^g = U_H^{c^nh} \times U_K^{c^nk} \le (C^p\ltimes H) \times K$$ because $C^p\ltimes H$ and $K$ are characteristic subgroups of~$(C^p\ltimes H)\times K$. Since all subgroups of~\hbox{$C\ltimes H$} are pronormal, there exists $u_1 \in \langle U_H, U_H^g \rangle \leq C^p\ltimes H$ such that $U_H^{u_1}=U_H^g$. On the other hand, by Theorem \ref{autpsl22p}, $U_K$ is a pronormal subgroup of $C\ltimes K$, so there exists $u_2 \in \langle U_K, U_K^g \rangle\leq K$ such that $U_H^{u_2}=U_H^g$. Since $u_1u_2\in\langle U,U^g\rangle$ and $U^{u_1u_2}=U_H^{u_1}\times U_{K}^{u_2}=U_H^g\times U_K^g=U^g$, the claim is proved.
 
 \medskip

\noindent {\it Step \textnormal(b\textnormal)}\quad Suppose that $|UC^pH/C^pH|$ is not a power of $2$ (so $UC^pH/C^pH$ is pronormal in~$G/C^pH$), and that~$U$ cannot be written as $U_{H}\times U_K$, where $U_H\leq C^p\ltimes H$ and $U_K\leq K$. We claim that $U$ is pronormal in~$G$.\\
\indent Let $g\in G$. Since \hbox{$\pi(S)\cap\pi(K)=\emptyset$,} Corollary 3.16 of \cite{ultimopaper} shows that we may first assume that $S=\{1\}$, and then that $L=C\times H$, where $H$ is a Dedekind $2$-group and $|C|=p$. Consequently, we may further assume that $U\cap H=\{1\}$, $H$ is an elementary abelian \hbox{$2$-group}, and $g\in C\ltimes K$; write $g=ck$, where $c\in C$ and $k\in K$. In particular, $U$ is finite, so, replacing $G$ by~\hbox{$\langle U,K,g,C\rangle$,} we may as well assume that $G$ is finite. We need to divide the proof according to the possibilities for~$U$. 

Suppose first that $p$ does not divide $|U|$, so $U\leq H\times K$. Since $U\cap H=\{1\}$, the hypotheses of~(b) yield that $U$ is dihedral (the only other possibility being that $U$ properly contains a~Sy\-low~\hbox{$2$-sub}\-group of $K$). Thus, $U=Y\ltimes N$, where $Y$ has order $2$ and $N$ is of odd order. As we have seen several times now (for example in the proof of point (b) of Theorem~\ref{plusautpsl22}), it is possible to assume that $N=N^g$. Write $Y=\langle xy\rangle$, where $x\in H$ and $y\in K$ are elements of order~$2$. Then $U^g=\langle xy^g\rangle\ltimes N$. Clearly, $Y^g\leq N_{H\times K}(N)=H\times N_K(N)$, so $\langle U,U^g\rangle\leq Y\ltimes N_K(N)$, and hence~$U$ and $U^g$ are (obviously) conjugate in~$\langle U,U^g\rangle$.

Now, suppose that $p$ divides $|U|$, so we may assume $C\leq U$. Write $U=C\ltimes N$, where $N\leq H\times K$. If $N$ is a pronormal subgroup of $G$, then $U$ is pronormal as well, being the product of two pronormal subgroups one of which normalizes the other. Thus, we may assume $N$ is not pronormal in $G$, and hence (by what we have already proved) that $N$ is an elementary abelian $2$-group. For a similar reason, we may further assume that $N$ is not normalized by $g$.


By hypothesis, $UH/H$ is pronormal in $G/H$, so there is $x\in\langle U,U^g\rangle$ with $U^xH=U^gH$. Thus,
$\langle U^x,U^g\rangle\leq \langle U,U^g\rangle$, and in order to prove that $U$ is conjugate to $U^g$ in $\langle U,U^g\rangle$, we may prove that $U^x$ is conjugate to $U^g$ in $\langle U^x,U^g\rangle$, so without loss of generality, we may assume $U^gH=UH$. Thus, $N^gH=NH$ and also $N^kH=NH$ because $c\in C\leq U$. 

Assume $k$ has odd order, so $NH/H$ is a Sylow $2$-subgroup of $G/H$. As we remarked above, $N^k=N^g\neq N$, so $\langle N,N^g\rangle\geq K$, and hence $\langle U,U^g\rangle\geq CK\ni g$. Thus, $k$ is an involution. But then $\langle k\rangle NH$ is a $2$-subgroup of $H\times K$ and so it is abelian. Consequently, $g$ normalizes $N$, and we have a contradiction. Point (b) is proved.

 \medskip

\noindent {\it Step \textnormal(c\textnormal)}\quad If $g$ is any involution of $K$, then $(C^p\ltimes H)\times \langle g\rangle$ is contained in the pro-norm of $G$.\\
\indent Let $U$ be a subgroup of $G$, and let $x\in (C^p\ltimes H)\times\langle g\rangle$. We need to show that $U$ and $U^x$ are conjugate in $\langle U,U^x\rangle$. To this aim we may assume $U$ is not pronormal, so by Steps (a) and~(b), $UC^pH/C^pH$ is a~\hbox{$2$-sub}\-group of~$G/C^pH$. Thus, $\langle g,U\rangle\leq (C^p\ltimes H)\times K$ and the statement follows from Step (c) of the proof of Theorem \ref{plusautpsl22}.

\medskip

By combining Steps (a), (b), and (c), the statement is proved.
\end{proof}

\medskip

In order to establish a general structure theorem that also encompasses infinite groups, we cannot rely exclusively on locally graded groups, as we did in the previous results. This difficulty arises, for instance, from the fact that although we know that infinite locally finite simple groups with our property do not exist (see Corollary~\ref{noinfinitesimple}), we are unable to prove an analogous result for infinite locally graded simple groups. To overcome this and related issues, we introduce the following definition: a group $G$ is said to be \emph{locally c-graded} if, in every finitely generated subgroup of $G$, the insoluble chief factors are locally finite. Clearly, every local\-ly~\hbox{c-graded} group is locally graded, and every locally (soluble-by-finite) group is local\-ly~\hbox{c-graded.} Moreover, the class of local\-ly~\hbox{c-graded} groups is obviously closed under the formation of subgroups and homomorphic images.

\begin{theo}\label{locallycgraded}
Let $G$ be a locally c-graded insoluble group. Then the subgroups of $G$ are pronormal or contained in the pro-norm of $G$ if and only if one of the following possibilities holds\textnormal:
\begin{itemize}
    \item[\textnormal{(1)}] $G=H\times K$, where\textnormal: 
    \begin{itemize}
        \item $H$ is periodic and all its subgroups are pronormal,
        \item $\pi(H)\cap\pi(K)\subseteq\{2\}$,
        \item if $Z(K)\neq\{1\}$, then  the Sylow $2$-subgroups of $H$ are elementary abelian, and
        \item $K$ is isomorphic to one of the following groups\textnormal: $J_1$, $\operatorname{PSL}(2,2^p)$, with $p$ a prime, or $\operatorname{PSL}(2,q)$, $\operatorname{SL}(2,q)$, where\textnormal: 
        \begin{itemize}
            \item $q=5^m$ with $m=1$ or $m$ an odd prime, or
            \item $q>5$ and $q\equiv_83,5$, or
            \item $q=2^m$, where $m$ is a prime.
        \end{itemize}
    \end{itemize}

    \item[\textnormal{(2)}] $G=C\ltimes (H\times K)$, where\textnormal: 
    \begin{itemize}
        \item $K=\operatorname{PSL}(2,2^p)$ for a prime $p>3$,
        \item $[C,K]\neq\{1\}$ and $C$ has order $p^n$ for some integer $n\geq1$,
        \item $\pi(H)\cap\pi(K)\subseteq\{2\}$,
        \item $H,K\trianglelefteq G$, and
        \item $H$ is periodic, and all subgroups of $C\ltimes H$ are pronormal.
    \end{itemize} 
\end{itemize}
\end{theo}
\begin{proof}
The sufficiency of the conditions follows by combining all the results we have obtained so far. Let us now consider the necessity. First we assume $G$ is finite. Let $M$ be a minimal normal subgroup of $G$. Suppose first that $M$ is an elementary abelian $p$-group for some prime $p$. By induction, $G/M$ has the required structure. In particular, there exists a finite normal non-abelian simple subgroup $K/M$ of~$G/M$. Let $X$ be a cyclic $2$-subgroup of $K$ such that~$XM/M$ is an involution. We know that~$XM/M$ is not pronormal in $K/M$, so~$X$ is not pronormal in $G$, and hence $X$ is contained in the pro-norm of $G$. Consequently,~$X$ normalizes every subgroup of $M$. The arbitrariness of~$X$ yields that every subgroup of $M$ is normal in $K$, so $M\leq Z(K)$ because $K/M$ is simple. By Theorem \ref{classificationsimple}, the Schur multiplier of~$K/M$ has order at most $2$. Consequently, $K=S\times R$, where $R$ is an elementary abelian \hbox{$p$-group,} and either $S\simeq K/M$ or $S$ is a non-splitting central extension of the simple group~$K/M$; in particular,~\hbox{$M=Z(S)\times R$.} But now $S=K'$, so $Z(S)$ is a characteristic subgroup of $K$. Thus, either $Z(S)\neq\{1\}$ and so $M=Z(S)$, or~$S$ is a simple normal subgroup of $G$ and $R=M$. Note that $Z(S)$ has order at most $2$, so $Z(S)\leq Z(G)$.

In both cases, let $C=C_G(S)$. By Theorem \ref{classificationsimple}, Corollary \ref{corpgl}, Theorem \ref{autpsl22p}, and~The\-o\-rem~\ref{aut5^p}, either $SC=G$, or $S=\operatorname{PSL}(2,2^p)$ and $|G:SC|=p>3$. Suppose the former. Then $S\cap C=Z(S)$. Since $S$ contains subgroups that are not contained in the pro-norm of $S$ (and so in that of $G$), $C/Z(S)$ has only pronormal subgroups. Thus, $C/Z(S)=U/Z(S)\ltimes V/Z(S)$, where $\pi(U/Z(S))\cap\pi(V/Z(S))=\emptyset$, $2\not\in\pi(V/Z(S))$, and $U/Z(S)$ is a Dedekind group acting as power automorphisms on $V/Z(S)$. We claim that $C$ splits over $Z(S)$. To this aim we may certainly assume that $|Z(S)|=2$, and that $C/Z(S)$ is a Dedekind $2$-group. Let $a$ be any element of $C$ such that $1\neq a^2\in Z(S)$, and let $Q$ be a~Sy\-low~\hbox{$2$-sub}\-group of $S$; in particular, $Q=\langle x,y\rangle$ is a quaternion group of order $8$. Now, $\langle ax\rangle$ is a subgroup of order $2$ that is not pronormal (as conjugation by $y$ shows) and is not even contained in the pro-norm (because $\langle ay\rangle$ is not pronormalized by $ax$). Thus, if $C/Z(S)$ has an element of order $4$, then $G$ has a subgroup of the form $\mathcal C_4\times S$, which is not possible by Corollaries \ref{corsl2,q}, \ref{corsl5q}, and \ref{corsl2q}. Thus, $C/Z(S)$ is an elementary abelian $2$-group, and such is also $C$. Consequently, $G=R\times S$ for some subgroup~$R$, and the statement follows from a combination of our previous results.

Suppose now that $S=\operatorname{PSL}(2,2^p)$ and $|G:SC|=p>3$. By the previous paragraph, we can write $SC=S\times R$ for some  subgroup $R$. Clearly, $p$ does not divide the order of $S$, and we may write $G=L\ltimes (S\times R)$, where $L$ is cyclic of $p$-power order. Then the result follows from~The\-o\-rem~\ref{damigliorare}.

Finally, the case in which $M$ is a non-abelian simple group can be handled exactly as above, by letting $M$ play the role of $S$. This proves the statement when $G$ is finite.

\medskip

Now, assume $G$ is locally finite. Since $G$ is not soluble, it has an insoluble finite subgroup~$E$ (see Corollary \ref{boundforsoluble}). By using the result in the finite case, we have that the soluble residual $R_E$ of $E$ is such that $R_E/Z(R_E)$ is simple non-abelian and $Z(R_E)$ has order at most $2$. If $F$ is any finite subgroup of $G$ containing $E$, then $F$ is insoluble and $R_F\geq R_E$. Thus, by assuming that the order of~$E$ is large enough, we may actually assume that $R_E=R_F$ for every finite subgroup~$F$ of $G$ containing $E$. It follows in particular that $R_E$ is normal in $G$. Finally, precisely as we did in the second and third paragraphs of the proof, we can finish the proof in this case.

\medskip

Assume $G$ is arbitrary. We claim that $G$ is locally finite, which will completely prove the statement. Suppose not, and let $F$ be an infinite finitely generated subgroup of $G$. Since $G$ is insoluble, it follows from Corollary \ref{boundforsoluble} that $G$ must contain a finitely generated insoluble subgroup, so without loss of generality we may assume that $F$ is insoluble. Let $R$ be the soluble residual of~$F$. It follows again from Corollary \ref{boundforsoluble} that $F/R$ is soluble, and so even finitely presented by~Co\-rol\-la\-ry~\ref{cortorsionfree}. This implies that $R$ is finitely generated as an~\hbox{$F$-in}\-va\-riant subgroup, and so that~$R$ contains a maximal~\hbox{$F$-in}\-variant subgroup~$E$. Since~$G$ is locally~\hbox{c-grad}\-ed, $R/E$ is either soluble, or insoluble and locally finite. The former is not possible (otherwise $R=E$), so~$R/E$ is locally finite and not soluble. But for locally finite groups we know the statement holds, so~$R/E$ is actually a finite simple non-abelian group. Clearly, $R/E\times C_{F/E}(R/E)$ has finite index in $F/E$. By combining all the results we have so far proved for direct products of simple groups with arbitrary subgroups, we have that $C_{F/E}(R/E)$ is periodic, so $F/E$ is locally finite, and hence even finite (being finitely generated). Let $E_1$ be any proper $F$-invariant subgroup of $E$ of finite index. Then $E/E_1$ is soluble because $F/E$ is insoluble. Moreover, the finite case yields that the soluble residual of $F/E_1$ is $R/E_1$, and this implies that $|E:E_1|\leq2$. But~$E_1$ is finitely generated and $G$ is locally graded, so $E_1=\{1\}$ and $F$ is finite, a contradiction.~\end{proof}

\begin{rem}\label{remfinale}
{\rm We can actually extended Theorem \ref{locallycgraded} to  smallest class of groups $\mathfrak D$ containing the class of locally c-graded groups that is local and it is closed under formation of ascending, descending series, and arbitrary subdirect products. Since the argument is very technical, based on the very abstract theory of classes described in the first chapter of \cite{Ro72}, and does not really extend the previous result in a very significant way, we only sketch the proof of this fact.

Let $G$ be a finitely generated group whose subgroups are pronormal or contained in the pro-norm. We need to show that if $G\in\mathfrak D$, then $G$ is either soluble or (finite simple non-abelian)-by-abelian. This can be done by transfinite induction by assuming that $G$ satisfies one of the following properties:
\begin{itemize}
    \item $G$ has an ascendant/descendant series whose factors are either soluble or (finite simple non-abelian)-by-abelian.
    \item $G$ has a family of normal subgroups whose quotient is either soluble or (finite simple non-abelian)-by-abelian, and whose intersection is trivial.
\end{itemize}

In each of these cases, recalling Corollary \ref{boundforsoluble}, we easily prove that $G$ has a finite normal series whose factors are either soluble or finite simple non-abelian, which implies that $G$ is a locally $c$-group, and so either soluble or (finite simple non-abelian)-by-abelian by our~The\-o\-rem~\ref{locallycgraded}.

\smallskip

By Lemma \ref{PLG->LF}, the class of locally graded $c$-groups contains the class of periodic locally graded groups, so the class $\mathfrak D$ contains the class of groups introduced by \v Cernikov in \cite{cernikov}.}
\end{rem}

\begin{cor}\label{coramel+}
Let $G$ be a locally graded group. If every subgroup of $G$ is pronormal or contained in the pro-norm, then the set of periodic elements forms a \textnormal(locally finite\textnormal) subgroup.
\end{cor}
\begin{proof}
Let $x,y \in G$ be periodic elements. Clearly, we may assume $G=\langle x,y\rangle$. Let $R$ be the finite residual of $G$. If $G/R$ is soluble, then it is periodic by Corollary \ref{cortorsionfree}, and so finite. If~$G/R$ is not soluble, then, being residually finite, it is (locally) finite by Theorem \ref{locallycgraded} and~Re\-mark~\ref{remfinale}. In both cases, $R=\{1\}$ and $G$ is finite.
\end{proof}

\begin{cor}
Let $G$ be an infinite locally graded group whose subgroups are pronormal or contained in the pro-norm. If $G$ is simple, then $G$ is torsion-free.
\end{cor}
\begin{proof}
This follows at once from Corollaries \ref{noinfinitesimple} and \ref{coramel+}.
\end{proof}

\bigskip\bigskip

\begin{flushleft}
\rule{8cm}{0.4pt}\\
\end{flushleft}

{
\sloppy
\noindent
Ernesto Ingrosso, Marco Trombetti

\noindent 
Dipartimento di Matematica e Applicazioni ``Renato Caccioppoli''

\noindent
Università di Napoli Federico II

\noindent
Complesso Universitario Monte S. Angelo

\noindent
Via Cintia, Napoli (Italy)

\noindent
e-mail: ernesto.ingrosso2@unina.it; marco.trombetti@unina.it 

}

\bigskip

{
\sloppy
\noindent
Amel Zitouni

\noindent 
Department of Sciences

\noindent
Teacher Education College of Setif-Messaoud Zeghar

\noindent
 19600, Setif (Algeria)

\noindent
e-mail: a.zitouni@ens-setif.dz

}

\end{document}